\documentclass[12pt,leqno]{article}
\usepackage[francais]{babel}
\usepackage{amssymb}
\usepackage{amsmath,amsfonts,amssymb}
\numberwithin{equation}{section}

\def\a{\`a}
\def\e{\`e}
\def\i{\`\i}

\begin{document}





${}$ 


\vskip 2 cm  

\centerline{\bf  An Initial and  boundary value problem on a strip } 
\centerline{\bf  for a large class of quasilinear hyperbolic systems arising} 
\centerline{\bf  from an atmospheric model.} 

\medskip 

\medskip 
\bigskip 
\medskip 
\centerline{{\bf Steave C. SELVADURAY${}^{1}$}} 

\bigskip 
\medskip 

\centerline{ 
${}^{1}$ Department of Mathematics ``GIUSEPPE PEANO''}
\centerline{Universit\`a di Torino, Torino, Italia}
\centerline{steave\_selva@yahoo.it, steaveclient.selvaduray@unito.it }

\bigskip 
\medskip 
\bigskip 
\medskip

{\small {\bf Abstract.} 
In this paper well-posedness is proved for an initial and boundary value problem (IBVP) relative to a large class of quasilinear hyperbolic systems, in $p+q$ equations, on a strip, arising from a model of $H_2O$-phase transitions in the atmosphere. To obtain this result, first, we extensively study an IBVP for the generic linear transport equation on $S_{t_1 }  = (0,t_1 ) \times G$ with uniformly locally Lipschitz data and associated vector field in $L_{x_d }^{-} (S_{t_1 } \times \mathbb{R}^p )$ (this cone of 
$L_t^\infty  (0,t_1 ;L_{(x,y_{loc} )}^\infty  (G \times \mathbb{R}^p ))^d 
$
is not cointained in $W^{1,\infty }_{\left( {t,x,y_{loc} } \right)}  \left( {S_{t_1}  \times \mathbb{R}^p} \right)^d$), that involves parametric vector functions  in $ L^\infty  \left( {0,t_1;W^{1,\infty } \left( G  \right)} \right)^p$,  by the method of characteristics. We obtain that the solution belongs to $W_{(t,x,y_{loc} )}^{1,\infty } (S_{t_1 }  \times \mathbb{R}^p )$ and interesting estimates about it.
 
Afterwards, using fixed point arguments, we establish the local existence, in time, and uniqueness of the solution in $W^{1,\infty } \left( {S_{t^* } } \right)^{p+q}$ for our class of quasilinear hyperbolic systems. Finally, we apply this result to study an IBVP for the hyperbolic part of an atmospheric model on the transition of water in the three states, introduced in \cite{[SF]}, such that rain and ice fall from it.

}
\bigskip 
\medskip

\bigskip 
\medskip
 

\bigskip 

\vskip 0.8 cm 

\section{- Introduction.}
In \cite{[SF]} we introduced and studied a model of motion of the air and the phase transitions of water in the atmosphere. The purpose of this research was to provide a detailed mathematical description about the phenomena which occur in the atmosphere such as wind or cloud formation and to show
its consistency. In order to obtain some mathematical results, we  introduced  several simplifications in our model. One of these assumed that the velocities of gas, water droplets and ice crystals are tangent to the boundary of a fixed spatial domain; therefore, we excluded, from our study, rain and snowfall. Hence, in the paper that follows,  we  include these phenomena limiting the study of the model to its hyperbolic part.

More precisely, the aim of the present paper is to establish a theorem of well-posedness, in  $W^{1,\infty }$, of an initial and  boundary value problem (IBVP) for a large class of quasilinear hyperbolic systems defined on a strip. There is an important application of this theorem, as we will see in the last section,  to an IBVP arising from a model of water phase transitions in the atmosphere given in \cite{[SF]}. Indeed, in the last section, we will obtain a result of well-posedness for the hyperbolic part of the atmospheric model defined in  a strip for which rain and ice fall from the strip   

Of course, to study this class of quasilinear hyperbolic systems, we use a linearization procedure. Therefore we need to have a result on the IBVP for the linear transport equation defined on a strip that can be applied to our quasilinear hyperbolic system.  The IBVP for the linear transport equation was studied in \cite{[B1]} by C. Bardos assuming time independent and Lipschitz vector field and using the method of characteristics and the semi-group theory. Hence, we can not use this result because the vector field of our system is depending on time.

Instead, in the paper \cite{[BO1]}, F. Boyer gives the trace theorems for the weak solutions of the linear transport equation on a regular bounded domain $\Omega$ and solves the corresponding IBVP, assuming in particular the vector field in $L^1 (0,T;W^{1,p} (\Omega))$ with null divergence and its trace in $L^p ((0,T) \times \partial \Omega )$ ($p>1$). Furthermore, O. Besson and J. Pousin study in \cite{[BP]} an IBVP for the linear transport equation on a bounded set assuming a $L^\infty-$vector field with $L^\infty-$divergence. They use a particular functional setting of an anisotropic Sobolev space and moreover the IBVP is reformulated by time-space least squares. Afterwards, G. Crippa, C. Donadello and L. V. Spinolo, establish in \cite{[CDS]}  well-posedness for continuity equations with bounded total variation coefficients; in particular, the vector field is bounded with bounded divergence and it is in $L_{loc}^1 \left( {0,T;BV(\Omega )} \right)$.      Unfortunately, in \cite{[BO1]}, \cite{[BP]} and \cite{[CDS]}, the domain is bounded and the vector field is also less regular to obtain estimates about the gradient of the solution for the linear transport equation, therefore we can not use these results  to study our IBVP for the quasilinear hyperbolic systems. Hence, using  a more regular vector field  depending also on time, the definition of generalized solution given in \cite{[M1]} and the method of characteristics, we obtain a  result on well-posedness in $W^{1,\infty }$ of the IBVP for the linear transport equation with parametric vector and source functions on a strip. Moreover,  from this result we deduce some useful estimates that play a vital role to study the analogous IBVP for our class of quasilinear hyperbolic systems. 

Let us say something about the sections of this paper. In section 2, we define some notations and functional spaces that will be used in that follows. In section 3, we introduce an IBVP for a large class of quasilinear hyperbolic systems and we transform this problem into system of integral equations of Volterra type using an extension of Cinquini-Cibrario method of characteristics. Therefore a generalized solution of our IBVP will be a solution of this system of integral equations. Hence, we give the statement of the main theorem regarding the well-posedness for the given IBVP about generalized solutions. In section 4, we introduce an analogous IBVP for a linear transport equation with vector parametric and source functions. In section 5, we study in detail the flow associated to the vector field relative to the linear transport equation, obtaining results about the regularity of the flow and the initial time of existence of the flow. From these results, we are able, in section 6, to prove the well-posedness of the IBVP for the linear transport equation that includes a parametric vector function. 

In section 7, after having obtained a preliminar result of existence and uniqueness  about the semilinear part of our quasilinear hyperbolic system, finally, we prove the main theorem of this paper. 

The paper ends (section 8) with the application of the main theorem to an IBVP for the hyperbolic part of an atmospheric model.   

The author would like to thank Prof. Hisao Fujita Yashima and  Prof. Davide Ascoli for some discussions on this subject. 
\medskip
\medskip
\medskip
\medskip
\section{- Some useful notations and functional spaces.}

In this section we define some notations and functional spaces that will be used later. First of all, let $G$ be an open subset of $\mathbb{R}^d$ defined as follows
\begin{equation}\label{ex(4-1-def-G)}
G=\mathbb{R}^{d-1} \times (0,1). 
\end{equation}
The canonical frame, in the linear space $\mathbb{R}^d$, is denoted by $\left\{ {e_j |j = 1,..,d} \right\}$. Moreover, we introduce  
\begin{equation}\label{ex(4-1-spazi-1)}
S_{t_1} =(0,t_1) \times G, \quad S'_{t_1} = (0,t_1)\times \mathbb{R}^{d-1}, 
\end{equation}
where $t_1>0.$ The generic point of $S_{t_1}$ is denoted by $(t,x)=(t,x',x_d)$. 

The maximum and minimum value of the real numbers $\alpha$ and $\beta$ is indicated by $\alpha \vee \beta$ and $\alpha \wedge \beta$ respectively; furthermore we define
\begin{equation}\label{ex(4-1-spazi-2)}
I\left( {\alpha ,\beta } \right) = \left( {\alpha  \wedge \beta } , {\alpha  \vee \beta }  \right) \quad \forall \alpha, \beta \in \mathbb{R}.
\end{equation}
Now, we introduce the following functional spaces
\begin{equation}\label{ex(4-1-spazi-4-1)}
L_t^r \big( {t_2,t_3;L_{(x,y_{loc} )}^{\infty } \left( {{G}  \times \mathbb{R}^p} \right)} \big) = 
\mathop  \cap \limits_{A > 0} L^r \left( {t_2,t_3;L_{}^{\infty } \left( {{G} \times (-A,A)^p } \right)} \right)
,\end{equation}
\begin{equation}\label{ex(4-1-spazi-4-5)}
L_t^r \big( {t_2,t_3;W_{(x,y_{loc} )}^{1,\infty } \left( {{G}  \times \mathbb{R}^p} \right)} \big) = 
\mathop  \cap \limits_{A > 0} L^r \left( {t_2,t_3;W_{}^{1,\infty } \left( {{G} \times (-A,A)^p } \right)} \right),
\end{equation}
where $0 \le t_2 \le t_3 \le t_1$, $1 \le r \le \infty$ and $p>0$ is an integer. Moreover, assuming $A>0$ and $0 \le a<b\le 1$,  we say that a measurable function $g:S_{t_1 }  \times \mathbb{R}^p \to \mathbb{R}
$ belongs to $L_{x_d }^r \big( {a,b;W_{(t,x',y )}^{1,\infty } \left( {S'_{t_1} \times (-A,A)^p } \right)} \big)$ if and only if
\begin{equation}\label{ex(4-1-integrale-uno-4-10)}
\left\| g \right\|^r_{L_{x_d }^r ( {a,b;W_{(t,x',y )}^{1,\infty } ( S'_{t_1} \times (-A,A)^p } ) )}  = \int\limits_a^b {\left\| {g( \cdot ,z, \cdot )} \right\|^r_{W^{1,\infty } (S'_{t_1}  \times ( - A,A)^p)} dz < \infty }; 
\end{equation}
therefore we can define
\begin{equation}\label{ex(4-1-spazi-4-11)}
L_{x_d }^r \big( {a,b;W_{(t,x',y_{loc} )}^{1,\infty } \left( {S'_{t_1}  \times \mathbb{R}^p} \right)} \big) =\mathop  \cap \limits_{A > 0} L_{x_d }^r \big( {a,b;W_{(t,x',y)}^{1,\infty } \left( {S'_{t_1} \times (-A,A)^p } \right)} \big).
\end{equation}
Furthermore the product of $k$ copies, for example, of $L_t^r \big( {t_2,t_3;W_{(x,y_{loc} )}^{1,\infty } \left( {{G}  \times \mathbb{R}^p} \right)} \big)$ will be denoted as follows 
\begin{equation}\label{ex(4-1-spazi-4-11-k)}
L_t^r \big( {t_2,t_3;W_{(x,y_{loc} )}^{1,\infty } \left( {{G}  \times \mathbb{R}^p} \right)} \big)^k ,
\end{equation}
where $k>0$ is an integer.

To avoid long estimates, we must introduce shorter symbols to indicate some norms. For example, we often consider a (vector) function $h:S_{t_1 }  \to \mathbb{R}^p$ (or $ \overline h $, or $h^{(l)}$, etc) with a given regularity and, in short, we put
\begin{equation}\label{ex(3-1-norma-h)}
\left\| h \right\|_\infty=\left\| h \right\|_{L^\infty  ( S_{t_1}  )}    ,\quad \left\| {\partial _{x_j } h} \right\|_\infty = \left\| {\frac{{\partial h}}{{\partial x_j }}
} \right\|_{L^\infty  ( S_{t_1} )}     ,\quad \left\| {D_{\left( {t,x} \right)} h} \right\|_\infty=\left\| {D_{\left( {t,x} \right)} h} \right\|_{L^\infty  ( S_{t_1} )}  ,  
\end{equation}
where $j=1,...,d$ and $D_{\left( {t,x} \right)} h$ is the jacobian matrix for $h$. Furthermore, if $\alpha$ and $\overline \alpha$ (vector or matrix functions on $S_{t_1}$)  have bounded components, then  we put
\begin{equation}\label{ex(3-1-norma-h-bis)}
\Lambda \left( {\alpha,\overline \alpha } \right) = 1 + \left\| \alpha \right\|_\infty   \vee \left\| {\overline \alpha } \right\|_\infty.
\end{equation}
In many situations, we consider a function $g \in L_t^r \big( {t_2,t_3;W_{(x,y_{loc} )}^{1,\infty } \left( {{G}  \times \mathbb{R}^p} \right)} \big)$,  two bounded functions $h,\overline h :S_{t_1 }  \to \mathbb{R}^p$ and we must study the following term 
\begin{equation}\label{ex(3-1-norma-a)}
I=\| \mathop {\rm ess \sup }\limits_{x \in G} g(\cdot,x,h(\cdot,x)) \vee \mathop {\rm ess \sup }\limits_{x \in G} g(\cdot,x,\overline h (\cdot,x))\|_{L^r (t_2 ,t_3 )}, 
\end{equation}
where $\rm ess \sup$ is the essential supremum.
It is obvious that we have
\begin{equation}\label{ex(3-1-norma-a)}
I \le \left\| g \right\|_{L^r (t_2 ,t_3 ;L^\infty  (G \times ( - \left\| h \right\|_\infty   \vee \left\| {\overline h } \right\|_\infty  ,\left\| h \right\|_\infty   \vee \left\| {\overline h } \right\|_\infty  )^p))}; 
\end{equation}
as well, for convenience, we introduce the following abbreviated notation
\begin{equation}\label{ex(3-1-norma-a1)}
 \left\| g \right\|_{L^r (t_2 ,t_3 ;h,\overline h )}=\left\| g \right\|_{L^r (t_2 ,t_3 ;L^\infty  (G \times ( - \left\| h \right\|_\infty   \vee \left\| {\overline h } \right\|_\infty  ,\left\| h \right\|_\infty   \vee \left\| {\overline h } \right\|_\infty  )^p))}, \quad 1 \le r \le \infty.
\end{equation}
In a similar way, we define
\begin{equation}\label{ex(3-1-norma-a2)}
\left\|  \cdot  \right\|_{L^\infty  (G;h,\overline h )}=\left\|  \cdot  \right\|_{L^\infty  (G \times ( - \left\| h \right\|_\infty   \vee \left\| {\overline h } \right\|_\infty  ,\left\| h \right\|_\infty   \vee \left\| {\overline h } \right\|_\infty  )^p)}   ,  
\end{equation}
\begin{equation}\label{ex(3-1-norma-a3)}
\left\|  \cdot  \right\|_{L^\infty  (S'_{t_1};h,\overline h )}=\left\|  \cdot  \right\|_{L^\infty  (S'_{t_1} \times ( - \left\| h \right\|_\infty   \vee \left\| {\overline h } \right\|_\infty  ,\left\| h \right\|_\infty   \vee \left\| {\overline h } \right\|_\infty  )^p)}   , 
\end{equation}
\begin{equation}\label{ex(3-1-norma-a4)}
 \left\| \cdot \right\|_{L^r_{x_d} (a ,b ;h,\overline h )}=\left\| \cdot  \right\|_{L^r_{x_d} (a ,b ;L^\infty  (S'_{t_1} \times ( - \left\| h \right\|_\infty   \vee \left\| {\overline h } \right\|_\infty  ,\left\| h \right\|_\infty   \vee \left\| {\overline h } \right\|_\infty  )^p))}.
\end{equation}
Of course, if $h= \overline h$, we can simplify the notations assuming  
\begin{equation}\label{ex(3-1-norma-a5)}
\left\| \cdot \right\|_{L^r (t_2 ,t_3 ;h )}=\left\| \cdot \right\|_{L^r (t_2 ,t_3 ;h,h )}, \quad \mbox{etc}.
\end{equation}
Moreover, assuming $K$ a subset of $\overline S_{t_1}$, we denote the family of uniformly locally lipschitz functions on $K$ by $Lip^{unif}_{loc}(K)$, i.e.  
\begin{equation}\label{ex(3-1-def-loc-lip)}
Lip^{unif}_{loc}(K)=\left\{ {\phi :K \to \mathbb{R}|\quad \exists M > 0:\mbox{ } \forall (t,x) \in K  } \right.\mbox{ } \exists \sigma  > 0:
\end{equation}
$$
\left. {\forall (t',x') \in K,\mbox{ } \left| {t - t'} \right| + \left| {x - x'} \right| \le \sigma \mbox{ } \Rightarrow \mbox{ }\left| {\phi (t,x) - \phi (t',x')} \right| \le M\left[ {\left| {t - t'} \right| + \left| {x - x'} \right|} \right]} \right\}.
$$

Finally, it is important in what follows to consider,  in 
$L_t^\infty  \big(0,t_1 ;L_{(x,y_{loc} )}^\infty  (G \times \mathbb{R}^p )\big )^d$, the following cone 
\begin{equation}\label{ex(4-1-spazi-4-12)}
L_{x_d }^{-} (S_{t_1 } \times \mathbb{R}^p ) = \Big\{b:S_{t_1 } \times \mathbb{R}^p  \to \mathbb{R}^d |\quad b(t,x,y)=\big(  b'(t,x,y),b_d( t,x) \big)  \mbox{ a.e. } (t,x,y) \in S_{t_1 } \times \mathbb{R}^p,
\end{equation}
$$
 b\in 
 L_t^\infty  \big(0,t_1 ;L_{(x,y_{loc} )}^\infty  (G \times \mathbb{R}^p )\big )^d \cap L_t^1 \big( {0,t_1;W_{(x,y_{loc} )}^{1,\infty } \left( {{G}  \times \mathbb{R}^p} \right)} \big)^d, \quad b_d \in L_{x_d }^1 \big( {0,1;W_{(t,x' )}^{1,\infty } \left( {S'_{t_1}  } \right)} \big) ,
$$
$$
\exists B_d >0: b_d(t,x) \le -B_d \quad \mbox{ a.e. } t \in (0,t_1),\quad  \forall x \in G \Big\} .
$$
{\bf Remark 2.0.0.} We observe that if $b \in L_{x_d }^{-} (S_{t_1 } \times \mathbb{R}^p )$ then $b_d$ can not belong to $W^{1,\infty } \left( {S_{t_1}  } \right)$. Indeed, assuming, for example $d=1$, we have
$$
\left( {x_1  - t - 1} \right)^{1/3}  - 1 \in L_{x_1 }^{-} (S_{t_1 } \times \mathbb{R}^p ) \backslash 
 W^{1,\infty } \left( {S_{t_1}  } \right).
$$ 

\medskip
\medskip
\medskip
\medskip
\section{- Position of the problem and main theorem.}

In this paper we study the IBVP for a large class of quasilinear hyperbolic systems in a strip, that generalizes the hyperbolic part, without integral terms, of a model of water phase transitions in the atmosphere given in \cite{[SF]} (see also Section 8). More precisely, we study,  in the unknown functions $y=(y_1,..,y_p)$, $w=(w_1,..,w_q)$, the following  IBVP
\begin{equation}\label{ex(1-1-IBVP-eq-y)}
\partial _t y_i (t,x)  + v_i \left( {t,x} \right) \cdot \nabla _x y_i(t,x)  = f_i \left( {t,x,y(t,x),w(t,x)} \right), \quad (t,x) \in S_{t_1},
\end{equation} 
\begin{equation}\label{ex(1-1-IBVP-eq-w)}
\partial _t w_j(t,x)  + u'_j \left( {t,x,y} \right) \cdot \nabla _{x'} w_j(t,x)  + u_{jd} \left( {t,x} \right)\partial _{x_d } w_j(t,x)  = 
\end{equation}
$$
=g_j \left( {t,x,y(t,x),w(t,x)} \right),\quad (t,x) \in S_{t_1},
$$
\begin{equation}\label{ex(1-1-dato-iniziale-y)}
y_i \left( {0,x} \right) = y_{0i} \left( x \right), \quad x \in G, 
\end{equation}
\begin{equation}\label{ex(1-1-dato-frontiera-w)}
w_j \left( {t,x} \right) = w_j^* \left( {t,x} \right), \quad (t,x) \in \Gamma_-,
\end{equation}
where $i=1,...,p$, $j=1,...,q$, 
\begin{equation}\label{ex(1-1-def-domain)}
\Gamma_- = \big( \left\{ 0 \right\} \times G \big)  \cup \big( \left[ {0,t_1 } \right] \times \mathbb{R}^{d-1} \times \left\{ 1 \right\} \big),
\end{equation}
\begin{equation}\label{ex(1-1-def-operatori)}
\partial_t= \frac{\partial }{{\partial t}}, \quad \partial_{x_k }= \frac{\partial }{{\partial x_k}}, \quad \nabla_x=\big(\partial x_1,...,\partial x_d \big), \quad \nabla_{x'}=\big(\partial x_1,...,\partial x_{d-1} \big),\quad k=1,...,d.
\end{equation}
Moreover  $f_i ,g_j :S_{t_1 }  \times \mathbb{R}^p  \times \mathbb{R}^q  \to \mathbb{R}$ are given sources for the corresponding transport equations, $v_i :S_{t_1 }    \to \mathbb{R}^d$ and $u_j=(u_{j1},u_{j2},...,u_{jd-1},u_{jd})=(u'_j,u_{jd}):S_{t_1 }  \times \mathbb{R}^p   \to \mathbb{R}^d$ are vector fields, $y_{0i}: G \to \mathbb{R}$,  $w^*_j: \Gamma_- \to \mathbb{R}$ are given functions and $i=1,..,p$, $j=1,..,q$. 

Afterwards, we can associate to the IBVP \eqref{ex(1-1-IBVP-eq-y)}-\eqref{ex(1-1-dato-frontiera-w)}, by the method of characteristics, the following system of integral equations 
\begin{equation}\label{ex(1-1-IS-eq-y)}
y_i (t,x) = y_{0i} ({\rm X}_{iY} (0;t,x)) + \int\limits_0^t {f_i (s,{\rm X}_{iY} (s;t,x),y(s,{\rm X}_{iY} (s;t,x),w(s,{\rm X}_{iY} (s;t,x)))ds}, 
\end{equation} 
\begin{equation}\label{ex(1-1-IS-eq-w)}
w_j (t,x) = w_j^* (\tau _{j - } (t,x,y),{\rm X}_{jW} (\tau _{j - } (t,x,y);t,x,y)) + 
\end{equation} 
$$
 + \int\limits_{\tau _{j - } (t,x,y)}^t {g_j (s,{\rm X}_{jW} (s;t,x,y),y(s,{\rm X}_{jW} (s;t,x,y)),w(s,{\rm X}_{jW} (s;t,x,y)))ds} ,
$$
where $i=1,..,p$, $j=1,..,q$, ${\rm X}_{iY}$ and ${\rm X}_{jW}$ are the fluxes associated to the vector fields $v_i$ and $w_j$ respectively and, therefore,  are defined as follows 
\begin{equation}\label{ex(1-1-IS-eq-X-Y)}
{\rm X}_{iY} (s;t,x) = x - \int\limits_s^t {v_i (r,{\rm X}_{iY} (r;t,x))dr}, \quad s \in [0,t], 
\end{equation} 
\begin{equation}\label{ex(1-1-IS-eq-X-W)}
{\rm X}_{jW} (s;t,x,y) = x - \int\limits_s^t {u_j (r,{\rm X}_{jW} (r;t,x,y),y(r,{\rm X}_{jW} (r;t,x,y)))dr},\quad s \in [\tau _{j - } (t,x,y),t], 
\end{equation} 
whereas ${\tau _{j - } (t,x,y)}$ is the minimal time of existence for the solution of \eqref{ex(1-1-IS-eq-X-W)}.

Therefore, we assume that a solution for the system 
\eqref{ex(1-1-IS-eq-y)} and \eqref{ex(1-1-IS-eq-w)}
(see Theorem 3.1) is a generalized solution for the IBVP \eqref{ex(1-1-IBVP-eq-y)}-\eqref{ex(1-1-dato-frontiera-w)}. This definition is consistent with the one given by A. Myshkis in \cite{[M1]}.

Hence, we assume the following conditions on the functions that appear in the IBVP \eqref{ex(1-1-IBVP-eq-y)}-\eqref{ex(1-1-dato-frontiera-w)}:
\begin{equation}\label{ex(1-10-f-g)}
f_i, g_j \in L^ \infty_{t } (0,t_1 ;L_{(x,y_{loc} ,w_{loc} )}^{\infty } (G  \times \mathbb{R}^p  \times \mathbb{R}^q )) \cap L^1_{t } (0,t_1 ;W_{(x,y_{loc} ,w_{loc} )}^{1,\infty } (G  \times \mathbb{R}^p  \times \mathbb{R}^q )),  
\end{equation} 
\begin{equation}\label{ex(1-10-condiz-vi)}
v_i \in L^\infty (0,t_1 ;L^{\infty } (G))^d\cap L^1 (0,t_1 ;W^{1,\infty } (G))^d,  
\end{equation}
$$
v_i(t,x',0) \cdot e_d = v_i(t,x',1) \cdot e_d=0 \quad \forall x' \in \mathbb{R}^{d-1},  \mbox{ a.e. } t \in (0,t_1),
$$
\begin{equation}\label{ex(1-10-condiz-uj)}
 u_j  \in L_{x_d }^ -  \left( {S_{t_1 }  \times \mathbb{R}^p } \right),
\end{equation}
\begin{equation}\label{ex(1-10-condiz-iniz)}
y_{0i}  \in W^{1,\infty } (G),\quad w_j^*  \in Lip^{unif}_{loc}(\Gamma _ -  ), \quad i=1,...,p, \quad j=1,...,q .
\end{equation}

Now, we are ready to give the main theorem of this paper
\medskip
\medskip

{\bf Theorem 3.1.}  \ {\it Assume that the hypotheses  \eqref{ex(1-10-f-g)}-\eqref{ex(1-10-condiz-iniz)}   are verified. Then there exists $0 < t^*  \le t_1$ such that the system of integral equations
\eqref{ex(1-1-IS-eq-y)} and \eqref{ex(1-1-IS-eq-w)}, where ${\rm X}_{iY}$ and ${\rm X}_{jW}$ are given by \eqref{ex(1-1-IS-eq-X-Y)} and \eqref{ex(1-1-IS-eq-X-W)}, admits one and only one solution $(y, w)\in$ 
$W^{1,\infty } \left( {S_{t^* } } \right)^{p+q}$.
The vector function $(y,w)$ is also said to be the generalized solution for IBVP \eqref{ex(1-1-IBVP-eq-y)}-\eqref{ex(1-1-dato-frontiera-w)}.

Moreover, for any sufficiently small t, the mapping $( y_0 ,w^* ,v_1 ,...,v_p , u_1,..., u_q,$ $f,g) \in
$ $W^{1,\infty } \left( {G } \right)^p  \times Lip^{unif}_{loc} \left( {\Gamma_- } \right)^q  \times \big( L^\infty (0,t ;L^{\infty } (G)) \cap L^1 \left( {{0,t };W^{1,\infty } \left( {G } \right)} \right) \big)^{dp} \times L_{x_d }^ -  \left( {S_{t }  \times \mathbb{R}^p } \right)^{q} \times 
\big [L^\infty_{t } (0,t ;L_{(x,y_{loc} ,w_{loc} )}^{\infty } (G  \times \mathbb{R}^p  \times \mathbb{R}^q ))  \cap 
L^1_{t } (0,t ;W_{(x,y_{loc} ,w_{loc} )}^{1,\infty } (G  \times \mathbb{R}^p  \times \mathbb{R}^q )) \big]^{p+q }$ $ \to \left( {y,w} \right) \in L^\infty  \left( {S_t} \right)^{p+q}$ is locally Lipschitz continuous.

}
\medskip
\medskip

To prove this theorem we must first deduce some results about the linear transport equations. 

\medskip
\medskip
\medskip
\medskip
\section{- An initial and boundary value problem for the linear transport equation with parametric vector function and source term.}

In this section we study the following linear transport equation 
\begin{equation}\label{ex(4-1-lin-transp-pde-mixed)}
\partial _t z\left( {t,x} \right) + b\left( {t,x,h(t,x)} \right) \cdot \nabla_x z\left( {t,x} \right)+ c\left( {t,x,h(t,x)} \right)z(t,x) = a\left( {t,x,h(t,x)} \right), \quad  \left( {t,x} \right) \in S_{t_1},
\end{equation} 
under the following condition
\begin{equation}\label{ex(4-1-condiz-boundary-value)}  
z \left( {t,x} \right) = z^* \left( {t,x} \right) \quad \left( {t,x} \right) \in \Gamma _ -  
, 
\end{equation} 
where $b=(b',b_d):S_{t_1}  \times \mathbb{R}^p  \to \mathbb{R}^d$ is a vector field with $b_d:S_{t_1}    \to \mathbb{R}$ ($d$ is an integer greater than zero) , $c:S_{t_1}\times \mathbb{R}^p \to \mathbb{R}$ and $a:S_{t_1}\times \mathbb{R}^p   \to \mathbb{R}$ (source) are  given functions, $h:S_{t_1}  \to \mathbb{R}^p$ is a parametric vector function, $\Gamma_-$ is the surface carrying data defined in \eqref{ex(1-1-def-domain)} and $z^* :\Gamma_- \to \mathbb{R}$ is a given function. The problem \eqref{ex(4-1-lin-transp-pde-mixed)}-\eqref{ex(4-1-condiz-boundary-value)} is called an initial and boundary value problem (IBVP) for the linear transport equation \eqref{ex(4-1-lin-transp-pde-mixed)}.

This IBVP will be studied through the analysis of the flow $\rm X$ of the vector field $b$ which satisfies
\begin{equation}\label{ex(4-1-equaz-caratteristica-mixed)}
\left\{ \begin{array}{l}
 \partial _t {\rm X}\left( {t;t_0,x_0,h} \right) = b\left( {t,{\rm X}\left( {t;t_0 ,x_0 ,h} \right),h\left( {t,{\rm X}\left( {t;t_0 ,x_0 ,h} \right)} \right)} \right)
 \\ 
 {\rm X}\left( {t_0;t_0,x_0,h} \right) = x_0 \quad (t_0,x_0) \in S_{t_1}.\\ 
 \end{array} \right.
\end{equation} 
In fact, we will see, under some hypotheses, that the solution $z$ at point $(t_0,x_0) \in S_{t_1}$ is determined by the knowledge of the flow ${\rm X}\left( {t;t_0,x_0,h} \right)$ which starts from $\tau _ -  \left( {t_0 ,x_0 ,h} \right)$ (initial time) and the value of $z^*$ in $\left( {\tau_- \left( {t_0,x_0,h} \right),{ X}\left( {\tau_- \left( {t_0,x_0,h} \right)}; t_0,x_0,h \right)} \right)
$ $\in \Gamma_{-}$ (see Section 6). The Cauchy problem \eqref{ex(4-1-equaz-caratteristica-mixed)} is called the characteristic problem related to the transport equation \eqref{ex(4-1-lin-transp-pde-mixed)}. Moreover, it is equivalent, in Caratheodory theory, to the following integral equation  
\begin{equation}\label{ex(4-1-equaz-integrale)}
{\rm X}\left( {t;t_0 ,x_0 ,h} \right) = x_0  - \int\limits_t^{t_0 } {b\left( {s,{\rm X}\left( {s;t_0 ,x_0 ,h} \right),h\left( {s,{\rm X}\left( {s;t_0 ,x_0 ,h} \right)} \right)
} \right)ds}. 
\end{equation}
\medskip
\medskip
\medskip
\medskip
\section{- Regularity of the flow ${\rm X}\left( {\cdot;t_0,x_0,h} \right)$ and the initial time $\tau _ -  \left( {t_0 ,x_0 ,h} \right)$.}

In this section  we give some useful lemmas about the flow ${\rm X}\left( {\cdot;t_0,x_0,h} \right)$ and the initial time $\tau _ -  \left( {t_0 ,x_0 ,h} \right)$.

{\bf Lemma 5.1.}   \ {\it Assume for the vector field $b$ the following regularity
\begin{equation}\label{ex(3-1-condiz-b-a)}
b \in L_t^1 \big( {0,t_1;W_{(x,y_{loc} )}^{1,\infty } \left( {{G}  \times \mathbb{R}^p} \right)} \big)^d 
\end{equation} 
and suppose that there exists a positive constant $B_d$ such that 
\begin{equation}\label{ex(3-1-condiz-bd)}
b_d(t,x) \le - B_d \quad \forall x \in G \quad \mbox{a.e. } t \in (0,t_1).
\end{equation} 
Then, for every $\left( {t_0 ,x_0,h } \right) \in S_{t_1} \times L^\infty  \left( {0,t_1;W^{1,\infty } \left( G  \right)} \right)^p
$, there exists
a unique maximal Caratheodory's solution ${\rm X}(\cdot;t_0,x_0,h):\left[ {\tau _ -  \left( {t_0 ,x_0 ,h} \right),\tau _ +  \left( {t_0 ,x_0 ,h} \right)} \right] \to \overline {G}$ of \eqref{ex(4-1-equaz-integrale)} with 
\begin{equation}\label{ex(4-1-punto-entrata)}
\left( {\tau _ -  \left( {t_0 ,x_0 ,h} \right),{\rm X}\left( {\tau _ -  \left( {t_0 ,x_0 ,h} \right);t_0 ,x_0 ,h} \right)} \right) \in \Gamma _ -. 
\end{equation}
Furthermore, the initial time $\tau _ -  \left( {t_0 ,x_0 ,h} \right)$ satisfies the following estimate
\begin{equation}\label{ex(4-1-durata-tempo-iniziale)}
0 \le \tau _ -  \left( {t_0 ,x_0 ,h} \right) \le t_0 - \rho,
\end{equation}
where $\rho$ is a number defined as follows 
\begin{equation}\label{ex(4-1-definizione-numero-d)}
0 < \rho \le t_0, \quad \left\| b \right\|_{L^1 (t_0  - \rho,t_0 ;h)}  < x_{0d}  \wedge \left( {1 - x_{0d} } \right).
\end{equation}

In the end, given $\left( {t_0 ,x_0 ,h} \right) \in S_{t_1}  \times L^\infty  \left( {0,t_1;W^{1,\infty } \left( G  \right)} \right)^p
$, the following statements hold
\begin {itemize}
{\bf \item i)} If $s, \overline s \in \big[ \tau _ -  \left( {t_0 ,x_0 ,h} \right),t_0 
\big]$ then
\begin{equation}\label{ex(3-1-norma-X-first)}
\left| {{\rm X}\left( {s;t_0 ,x_0 ,h} \right) - {\rm X}\left( {\overline s ;t_0 ,x_0 ,h} \right)} \right| \le \int\limits_{I\left( {s,\overline s } \right)} {\left\| {b\left( {r, \cdot } \right)} \right\|_{L^\infty  \left( {{G} ;{h  } } \right)} dr}. 
\end{equation}

If $\left( {\overline t _0 ,\overline x _0 ,\overline h } \right) \in S_{t_1}  \times L^\infty  \left( {0,t_1;W^{1,\infty } \left( G  \right)} \right)^p
$ and $s \in \big[  \tau _ -  \left( {t_0 ,x_0 ,h} \right) \vee \tau _ -  \left( {\overline t _0 ,\overline x _0 ,\overline h } \right),t_0  \wedge \overline t _0 \big]
$ then
\begin{equation}\label{ex(3-1-norma-X-second)}
\left| {{\rm X}\left( {s;t_0 ,x_0 ,h} \right) - {\rm X}\left( {s;\overline t _0 ,\overline x _0 ,\overline h } \right)} \right| \le 
\end{equation}
$$
\le C_1\Big( {{ \left\| {b} \right\|_{L^1 (I(t_0 ,\overline t _0 )
;h, \overline h)} + \left| {x_0  - \overline x _0 } \right| + \left\| {h - \overline h } \right\|_\infty  }  } \Big),
$$
where
\begin{equation}\label{ex(3-1-def-C-1)}
C_1=\left( {1 \vee \left\| {D_{\left( {x,y} \right)} b} \right\|_{L^1 ( 0,t_1;h,\overline h )} } \right)\exp \left( {\Lambda \left( {D _x h ,D _x \overline h  } \right)\left\| {D_{\left( {x,y} \right)} b} \right\|_{L^1 ( 0,t_1;h,\overline h )} } \right).
\end{equation}
$($$\Lambda$ is defined in \eqref{ex(3-1-norma-h-bis)}$)$.
{\bf \item ii)} Let   
\begin{equation}\label{ex(3-1-condiz-b-c)}
b \in 
L_t^\infty  \big(0,t_1 ;L_{(x,y_{loc} )}^\infty  (G \times \mathbb{R}^p )\big )^d,
\end{equation}
then  \eqref{ex(3-1-norma-X-first)} and \eqref{ex(3-1-norma-X-second)} can be replaced, respectively, by
\begin{equation}\label{ex(3-1-norma-X-first-bis)}
\left| {{\rm X}\left( {s;t_0 ,x_0 ,h} \right) - {\rm X}\left( {\overline s ;t_0 ,x_0 ,h} \right)} \right| \le 
\left\| b \right\|_{L^\infty  \left( {0,t_1;h } \right)} \left| {s - \overline s } \right|,
\end{equation}
\begin{equation}\label{ex(3-1-norma-X-second-bis)}
\left| {{\rm X}\left( {s;t_0 ,x_0 ,h} \right) - {\rm X}\left( {s;\overline t _0 ,\overline x _0 ,\overline h } \right)} \right| \le C_2\Big[ { \left| {t_0  - \overline t _0 } \right| + \left| {x_0  - \overline x _0 } \right| + \left\| {h - \overline h } \right\|_\infty  }  \Big],
\end{equation}
where
\begin{equation}\label{ex(3-1-def-C-2)}
C_2  = \left( {1 \vee \left\| {b} \right\|_{L^\infty  ( 0,t_1;h,\overline h )}  \vee \left\| {D_{\left( {x,y} \right)} b} \right\|_{L^1 (0,t_1;h,\overline h )} } \right)\times
\end{equation}
$$
\times \exp \left( {\Lambda \left( {D _x h,D _x \overline h  } \right)\left\| {D_{\left( {x,y} \right)} b} \right\|_{L^1 (0,t_1;h,\overline h )} } \right).
$$
\end {itemize}

}

\medskip

{\sc {\bf Proof}}. \

The integral equation \eqref{ex(4-1-equaz-integrale)} admits one and only one maximal solution ${\rm X} \in W_{loc}^{1,1} \left( {\tau _ -  \left( {t_0 ,x_0 ,h} \right),\tau _ +  \left( {t_0 ,x_0 ;h} \right)} \right)
$ (see \cite{[FV]}). Now, we show that ${\rm X}$ can be extended continuously on the closed interval $\left[ {\tau _ -  \left( {t_0 ,x_0 ;h} \right),\tau _ +  \left( {t_0 ,x_0 ;h} \right)} \right]$ and therefore it also satisfies \eqref{ex(4-1-equaz-integrale)} at the end points of this interval. Indeed, assuming $s_1$, $s_2$ such that $\tau _ -  \left( {t_0 ,x_0 ,h} \right) < s_1 \le s_2 < \tau _ + \left( {t_0 ,x_0 ;h} \right)$
, we have 
\begin{equation}\label{ex(3-2-stima-slz-1)}
\left| {{\rm X}\left( {s_2 ;t_0 ,x_0 ,h} \right) - {\rm X}\left( {s_1 ;t_0 ,x_0 ,h} \right)} \right| \le \int\limits_{s_1 }^{s_2} {\left\| {b\left( {r, \cdot } \right)} \right\|_{L^\infty  \left( {{G};h } \right)
}
 dr}, 
\end{equation}
then, thanks to a Cauchy criterion, we deduce that exists 
\begin{equation}\label{ex(3-2-limiti)}
 \mathop {\lim }\limits_{s \to \tau _ -  \left( {t_0 ,x_0 ,h} \right)^ +  } {\rm X}\left( {s;t_0 ,x_0 ,h} \right) \in \overline {G}. 
\end{equation} 
In a similar way we study the left limit of ${\rm X}\left( { s ;t_0 ,x_0 ,h} \right)$ as $s$ approaches $\tau _ +   \left( {t_0 ,x_0 ,h} \right)$. Therefore ${\rm X}\left( { \cdot ;t_0 ,x_0 ,h} \right) \in W^{1,1} \left[ {\tau _ -  \left( {t_0 ,x_0 ,h} \right),\tau _ +   \left( {t_0 ,x_0 ,h} \right)} \right]$. Taking into account \eqref{ex(3-1-condiz-bd)}, it is easy to get \eqref{ex(4-1-punto-entrata)}. To prove \eqref{ex(4-1-durata-tempo-iniziale)}, see, for example, \cite{[FV]} (Remark at page 5).     Furthermore, the estimate \eqref{ex(3-1-norma-X-first)} immediately follows from \eqref{ex(4-1-equaz-integrale)}.  

Now, given $\left( {t_0 ,x_0 ,h} \right),\left( {\overline t _0 ,\overline x _0 ,\overline h } \right) \in S_{t_1}  \times L^\infty  \left( {0,t_1;W^{1,\infty } \left( G  \right)} \right)^p
$, we consider the solutions ${\rm X}\left( { \cdot ;t_0 ,x_0 ,h} \right),{\rm X}\left( { \cdot ;\overline t _0 ,\overline x _0 ,\overline h } \right)$ and assume  $\tau _ -  \left( {t_0 ,x_0 ,h} \right) \vee \tau _ -  \left( {\overline t _0 ,\overline x _0 ,\overline h } \right) \le  t_0  \wedge \overline t _0$. The difference between the integral representations of the previous solutions (see \eqref{ex(4-1-equaz-integrale)}) gives, if $\tau _ -  \left( {t_0 ,x_0 ,h} \right) \vee \tau _ -  \left( {\overline t _0 ,\overline x _0 ,\overline h } \right) \le s \le  t_0  \wedge \overline t _0$, the following estimate 
\begin{equation}\label{ex(3-2-magg-xd)}
\left| {X_d \left( {s;t_0 ,x_0 ,h} \right) - X_d \left( {s;\overline t _0 ,\overline x _0 ;\overline h } \right)} \right| \le \left| {x_{0d}  - \overline x _{0d} } \right| + \int\limits_{I(t_0 ,\overline t _0 )} {{\left\| {b_d \left( {r, \cdot } \right)} \right\|_{L^\infty  \left( {G}  \right)} } dr +} 
 \end{equation}
$$
 + \int\limits_s^{t_0  \wedge \overline t _0 } {\left| {b_d \left( {r,X \left( {r;t_0 ,x_0 ,h} \right)} \right) - b_d \left( {r,X \left( {r;\overline t _0 ,\overline x _0 ;\overline h } \right)} \right)} \right|} dr,
$$
from which we immediately deduce
\begin{equation}\label{ex(3-2-magg-xd-10)}
\left| {X_d \left( {s;t_0 ,x_0 ,h} \right) - X_d \left( {s;\overline t _0 ,\overline x _0 ;\overline h } \right)} \right| \le \left| {x_{0d}  - \overline x _{0d} } \right| + \left\| {b_d } \right\|_{L^1 (I(t_0 ,\overline t _0 )
,L^\infty  \left( {G}  \right))}  + 
\end{equation}
$$
 + \int\limits_s^{t_0  \wedge \overline t _0 } {\left\| {\nabla _x b_d \left( {r, \cdot } \right)} \right\|_{L^\infty  \left( {G}  \right)} \left| {X \left( {r;t_0 ,x_0 ,h} \right) - X \left( {r;\overline t _0 ,\overline x _0 ;\overline h } \right)} \right|} dr.
$$
Analogously we obtain  
\begin{equation}\label{ex(3-2-magg-x'-20)}
\left| {X'\left( {s;t_0 ,x_0 ,h} \right) - X'\left( {s;\overline t _0 ,\overline x _0 ;\overline h } \right)} \right| \le \left| {x'_0  - \overline x '_0 } \right| + \left\| {b'} \right\|_{L^1 (I(t_0 ,\overline t _0 )
;h, \overline h)}  + 
\end{equation} 
$$
+ \left\| {D_{y} b'} \right\|_{L^1 ( {0,t_1} ;h,\overline h)} \left\| {h - \overline h } \right\|_{\infty}  + \int\limits_s^{t_0  \wedge \overline t _0 } \big( \left\| {D_{x} b'\left( {r, \cdot  } \right)} \right\|_{L^\infty  ({G};h, \overline h )}+\left\| {D_{y} b'\left( {r, \cdot  } \right)} \right\|_{L^\infty  ({G};h, \overline h )}\left\| {D _x h} \right\|_{\infty} \big) \times 
$$
$$
\quad \quad \quad \quad \quad \times \left| {X \left( {r;t_0 ,x_0 ,h} \right) - X \left( {r;\overline t _0 ,\overline x _0 ;\overline h } \right)} \right| dr.
$$
Adding \eqref{ex(3-2-magg-xd-10)}-\eqref{ex(3-2-magg-x'-20)} and using Gronwall's lemma, we deduce \eqref{ex(3-1-norma-X-second)}.

Of course, at this point, using again \eqref{ex(3-2-magg-x'-20)}, it is not so difficult to show {\bf ii)}. 

\  $ \square$ 
\medskip
\medskip

Now, assuming a more regular vector field $b$, we can obtain estimates on $\tau _ -  \left( {t_0 ,x_0 ,h} \right)$ and ${\rm X}\big( {\tau _ -  \left( {t_0 ,x_0 ,h} \right);t_0 ,x_0 ,h} \big)$.

{\bf Lemma 5.2.}   \ {\it Let us assume 
$\left( {t_0 ,x_0 ,h} \right) \in S_{t_1}  \times L^\infty  \left( {0,t_1;W^{1,\infty } \left( G  \right)} \right)^p
$. Then, 
if  
\begin{equation}\label{ex(3-1-condiz-bd-c)}
b \in L_{x_d }^{-} (S_{t_1} \times \mathbb{R}^p ) \quad \quad (\mbox{see \eqref{ex(4-1-spazi-4-12)}}),
\end{equation}
then there exists $\delta>0$ such that, under the conditions
$$\left( {\overline t_0,\overline x_0,\overline h } \right) \in S_{t_1} \times L^\infty  \left( {0,t_1;W^{1,\infty } \left( G  \right)} \right)^p,$$
$$
\left| {x_0  - \overline x _0 } \right|\le \delta, \quad \left| {x_0  - \overline x _0 } \right|\le \delta, \quad \left\| {h - \overline h } \right\|_\infty   \le \delta,
$$
we have
\begin{equation}\label{ex(3-1-norma-tau-first)}
\left| {\tau _ -  \left( {t_0 ,x_0 ,h} \right) - \tau _ -  \left( {\overline t_0 ,\overline x _0 ,\overline h } \right)} \right| + 
\end{equation}
$$
+\left| {{\rm X}\left( {\tau _ -  \left( {t_0 ,x_0 ,h} \right);t_0 ,x_0 ,h} \right) - {\rm X}\left( {\tau _ -  \left( {\overline t_0 ,\overline x _0 ,\overline h } \right);\overline t_0 ,\overline x _0 ,\overline h } \right)} \right| \le
$$
$$
\le C_3 \left[ {\left| {t_0  - \overline t _0 } \right| +\left| {x_0  - \overline x _0 } \right| + \left\| {h - \overline h } \right\|_\infty  } \right],
$$
where  
\begin{equation}\label{ex(3-1-def-C-3)}
C_3  = C\left( {1 + \left\| {b} \right\|_{L^\infty  (0,t_1;h,\overline h )}^2 +\left\| {D_{\left( {x,y} \right)} b} \right\|_{L^1 (0,t_1;h,\overline h )}^2 } \right)\times
\end{equation}
$$
\times \exp \left[ {C\left( {\left\| {b} \right\|_{L^\infty  (0,t_1;h,\overline h )}  + \Lambda \left( {D _x h,D _x \overline h } \right)\left\| {D_{\left( {x,y} \right)} b} \right\|_{L^1 (0,t_1;h,\overline h )} } \right)} \right]
$$
with a constant $C$ depending only on $B_d$, $\left\| {D _{\left( {t,x'} \right)} b_d } \right\|_{L_{x_d }^1 ( {0,1;L_{(t,x')}^\infty  (S_{t_1}' )} )} $.

Furthermore, for every $s \in \left( {\tau _ -  \left( {t_0 ,x_0 ,h} \right),t_0 } \right)
$, there exists $\delta>0$ such that, under the conditions
$$\left( {\overline t _0 ,\overline x _0 ,\overline h } \right) \in S_{t_1}  \times W^{1,\infty } \left( {S_{t_1} } \right)^p,\quad \overline s \in \left( {\tau _ -  \left( {\overline t_0 ,\overline x_0 , \overline h} \right),\overline t_0 } \right) ,  $$
$$
\left| {t_0  - \overline t _0 } \right| \le \delta, \quad \left| {x_0  - \overline x _0 } \right|\le \delta, \quad \left\| {h - \overline h } \right\|_{\infty}\le \delta,\quad \left| {s  - \overline s  } \right|  \le  \delta ,
$$
we have
\begin{equation}\label{ex(4-1-terza-stima-X)}
\left| {{\rm X}\left( {s;t_0 ,x_0 ,h} \right) - {\rm X}\left( {\overline s;\overline t _0 ,\overline x _0 ,\overline h } \right)} \right| \le C_4 \left[ { \left| {t_0  - \overline t _0 } \right|+ \left| {s  - \overline s } \right|+\left| {x_0  - \overline x _0 } \right| + \left\| {h - \overline h } \right\|_\infty   } \right],
\end{equation}
where
\begin{equation}\label{ex(4-1-def-C-5)}
C_4  = C\Big[ 1+ \left\| b \right\|_{L^\infty  (0,t_1;h,\overline h )} + \left( {1 \vee \left\| {b} \right\|_{L^\infty  ( 0,t_1;h,\overline h )}  \vee \left\| {D_{\left( {x,y} \right)} b} \right\|_{L^1 (0,t_1;h,\overline h )} } \right)\left\| {D_{\left( {x,y} \right)} b} \right\|_{L^1 (0,t_1;h,\overline h)}\times 
\end{equation}
$$
\times{\Lambda \left( {D _x h,D _x \overline h } \right)}  \Big]  \exp \left( {\Lambda \left( {D _x h,D _x \overline h } \right)\left\| {D_{\left( {x,y} \right)} b} \right\|_{L^1 (0,t_1;h,\overline h )} } \right),
$$
$($$C$ is a positive constant not depending on $b,h$, etc. $)$.
  
Hence, $\tau _ -   \in W_{(t,x,y_{loc} )}^{1,\infty } (S_{t_1 }  \times \mathbb{R}^p )
$, ${\rm X} \circ \tau _ -  \in W_{(t,x,y_{loc} )}^{1,\infty } (S_{t_1 }  \times \mathbb{R}^p )^d $.

}

\medskip

{\sc {\bf Proof}}. \ Now, assuming \eqref{ex(3-1-condiz-bd-c)}, we observe that ${\rm X}_d \left( { \cdot ;t_0 ,x_0 ,h} \right)
:\left[ {\tau _ -  \left( {t_0 ,x_0 ,h} \right),\tau _ +  \left( {t_0 ,x_0 ,h} \right)} \right]$ $ \to \left[ {{\rm X}_d \left( {\tau _ +  \left( {t_0 ,x_0 ,h} \right);t_0 ,x_0 ,h} \right),{\rm X}_d \left( {\tau _ -  \left( {t_0 ,x_0 ,h} \right);t_0 ,x_0 ,h} \right)} \right]$ is absolute continuous and strictly decreasing; therefore, the inverse function of
 ${\rm X}_d$, which we denote by ${\rm T}\left( { \cdot ;t_0 ,x_0 ,h} \right)$, is  continuous and  strictly decreasing. Hence, it is possible to apply the chain rule to the following function 
\begin{equation}\label{ex(3-30-def-id)}
id|_{\left[ {{\rm X}_d \left( {\tau _ +  \left( {t_0 ,x_0 ,h} \right);t_0 ,x_0 ,h} \right),{\rm X}_d \left( {\tau _ -  \left( {t_0 ,x_0 ,h} \right);t_0 ,x_0 ,h} \right)} \right]}  = {\rm X}_d  \circ {\rm T},
\end{equation} 
(see \cite{[LG]}). Therefore, we obtain:
\begin{equation}\label{ex(3-30-def-Td)}
\partial _{x_d } \Phi \left( {x_d } \right) = \frac{1}{{b_d \left( {\Phi \left( {x_d } \right),{\rm X}'\left( {\Phi \left( {x_d } \right);t_0 ,x_0 ,h} \right),x_d } \right)}}
\end{equation} 
for $x_d-$almost everywhere on $ \left[ {{\rm X}_d \left( {\tau _ +  \left( {t_0 ,x_0 ,h} \right);t_0 ,x_0 ,h} \right),{\rm X}_d \left( {\tau _ -  \left( {t_0 ,x_0 ,h} \right);t_0 ,x_0 ,h} \right)} \right]$, where
\begin{equation}\label{ex(3-30-def-PHI)}
\Phi \left( {\cdot } \right) = {\rm T}\left( {\cdot ;t_0 ,x_0 ,h} \right).
\end{equation}
Hence, we can associate to \eqref{ex(3-30-def-Td)} the following integral equation
\begin{equation}\label{ex(3-30-def-eq-int-T)}
\Phi\left( {x_d } \right) = t_0  + \int\limits_{x_{0d} }^{x_d } {\frac{{dw}}{{b_d \left( {\Phi \left( w \right),{\rm X}'\left( {\Phi \left( w \right);t_0 ,x_0 ,h} \right),w} \right)
 }}}. 
\end{equation}
for every $x_d  \in \left[ {{\rm X}_d \left( {\tau _ +  \left( {t_0 ,x_0 ,h} \right);t_0 ,x_0 ,h} \right),{\rm X}_d \left( {\tau _ -  \left( {t_0 ,x_0 ,h} \right);t_0 ,x_0 ,h} \right)} \right]$. 

Afterwards, studying the equation \eqref{ex(3-30-def-eq-int-T)} we deduce a key estimate to prove this lemma. Indeed, assuming $\left( {\overline t_0 ,\overline x_0 ,\overline h} \right) \in S_{t_1}  \times L^\infty  \left( {0,t_1;W^{1,\infty } \left( G  \right)} \right)^p
$, the condition \eqref{ex(3-1-condiz-bd-c)} 
and that exists $\beta$ such that $x_{0d}  \vee \overline x _{0d}  \le \beta  \le {\rm X}_d \left( {\tau _ -  \left( {t_0 ,x_0 ,h} \right);t_0 ,x_0 ,h} \right) \wedge {\rm X}_d \left( {\tau _ -  \left( {\overline t _0 ,\overline x _0 ,\overline h } \right);\overline t _0 ,\overline x _0 ,\overline h } \right)$, then, after having defined $\overline \Phi  \left(  \cdot  \right) = {\rm T}\left( { \cdot ;\overline t _0 ,\overline x _0 ,\overline h } \right)$ and assumed $x_d \in \left[ {x_{0d}  \vee \overline x _{0d} ,\beta } \right]$,  we obtain  
\begin{equation}\label{ex(3-30-stime-T)}
\left| {\Phi \left( {x_d } \right) - \overline \Phi \left( {x_d  } \right)} \right| \le 
\end{equation}
$$
\le \left| {t_0  - \overline t _0 } \right| +  B_d ^{ - 1} \left| {x_{0d}  - \overline x _{0d} } \right|  
+  B_d ^{ - 2} \int\limits_{x_{0d}  \vee \overline x _{0d} }^{x_d } \left\| {\nabla _{\left( {t,x'} \right)} b_d \left( { \cdot ,w} \right)} \right\|_{L^\infty  \left( {S'_{t_1} } \right)}\times 
$$
$$
 \times \left[ {\left| {\Phi \left( {w } \right)  - \overline \Phi \left( {w } \right)  } \right| + \left| {{\rm X}'\left( {\Phi \left( {w } \right) ;t_0 ,x_0 ,h} \right) - {\rm X}'\left( {\overline \Phi \left( {w } \right)  ;\overline t _0 ,\overline x _0 ,\overline h } \right)} \right|} \right] dw. 
$$
Now, for every $M>0$ there exists $\overline \rho>0$ such that if $\left| {x_0  - \overline x _0 } \right| \le M$, $\left\| {h - \overline h } \right\|_\infty   \le M$ and $\left| {t_0  - \overline t _0 } \right| \le \overline \rho$, then we have
\begin{equation}\label{ex(3-3-stime-tempi-iniziali-bis)}
\tau _\_ (t_0 ,x_0 ,h) \vee \tau _\_ (\overline t _0 ,\overline x _0 ,\overline h ) \le t_0  \wedge \overline t _0. 
\end{equation}
This result can be obtained by \eqref{ex(4-1-durata-tempo-iniziale)}. Therefore, assuming $\overline t_0 \le t_0$ and $ x_d \in \left[ {x_{0d}  \vee \overline x _{0d} ,\beta } \right]$, we deduce the following possibilities

1)$\tau _\_ (t_0 ,x_0 ,h)  \le \tau _\_ (\overline t _0 ,\overline x _0 ,\overline h ).$

   Hence, thanks to \eqref{ex(3-1-norma-X-first-bis)} and \eqref{ex(3-1-norma-X-second-bis)}, we obtain   
\begin{equation}\label{ex(3-3-magg-X'-first-0)}
\left| {{\rm X}'\left( {\Phi \left( w \right);t_0 ,x_0 ,h} \right) - {\rm X}'\left( {\overline \Phi  \left( w \right);t_0 ,\overline x _0 ,\overline h } \right)} \right| \le  
\end{equation}
$$
\le \left| {{\rm X}(\Phi (w);t_0 ,x_0 ,h) - {\rm X}(\overline \Phi  (w);t_0 ,x_0 ,h)} \right| + \left| {{\rm X}(\overline \Phi  (w);t_0 ,x_0 ,h) - {\rm X}(\overline \Phi  (w);\overline t _0 ,\overline x _0 ,\overline h )} \right| \le
$$
$$
\le \left\| b \right\|_{L^\infty  (0,t_1; h )} \left| {\Phi \left( w \right) - \overline \Phi  \left( w \right)} \right| + C_2 \left[ {\left| {t_0  - \overline t _0 } \right| +\left| {x_0  - \overline x _0 } \right| + \left\| {h - \overline h } \right\|_\infty  } \right],
$$
$\mbox{for all } w \in \left[ {x_{0d}  \vee \overline x _{0d} ,x_d } \right].
$
 
2)  $\tau _\_ (\overline t _0 ,\overline x _0 ,\overline h ) \le \tau _\_ (t_0 ,x_0 ,h).$

It is not restrictive to assume ${\overline \Phi  (w)} \le \tau _\_ (t_0 ,x_0 ,h)$ and $\overline t_0 \le \Phi  (w)$ with $w \in \left[ {x_{0d}  \vee \overline x _{0d} ,x_d } \right].$ therefore, thanks to \eqref{ex(3-1-norma-X-first-bis)} and \eqref{ex(3-1-norma-X-second-bis)}, we have the following estimate
\begin{equation}\label{ex(3-3-magg-X'-first)}
\left| {{\rm X}'\left( {\Phi \left( w \right);t_0 ,x_0 ,h} \right) - {\rm X}'\left( {\overline \Phi  \left( w \right);t_0 ,\overline x _0 ,\overline h } \right)} \right| \le  
\end{equation}
$$
\le \left| {{\rm X}(\Phi (w);t_0 ,x_0 ,h) - {\rm X}(\overline t_0;t_0 ,x_0 ,h)} \right| + \left| {{\rm X}(\overline t_0;t_0 ,x_0 ,h) - {\rm X}(\overline t_0;\overline t _0 ,\overline x _0 ,\overline h )} \right| \le
$$
$$
+\left| {{\rm X}(\overline t_0;t_0 ,x_0 ,h) - {\rm X}(\overline \Phi  (w);\overline t _0 ,\overline x _0 ,\overline h )} \right| \le \left\| b \right\|_{L^\infty  (0,t_1; h )}(\Phi (w) - \overline t _0 ) + 
$$
$$
+ C_2 \left[ {\left| {t_0  - \overline t _0 } \right| +\left| {x_0  - \overline x _0 } \right| + \left\| {h - \overline h } \right\|_\infty  } \right] + \left\| b \right\|_{L^\infty  (0,t_1;\overline  h )}( \overline t _0-\overline \Phi (w)  ) \le
$$
$$
\le \left\| b \right\|_{L^\infty  (0,t_1; h, \overline h )} \left| {\Phi \left( w \right) - \overline \Phi  \left( w \right)} \right| + C_2 \left[ {\left| {t_0  - \overline t _0 } \right| +\left| {x_0  - \overline x _0 } \right| + \left\| {h - \overline h } \right\|_\infty  } \right],
$$
$\mbox{for all } w \in \left[ {x_{0d}  \vee \overline x _{0d} ,x_d } \right].
$

In a similar way, we obtain the estimate \eqref{ex(3-3-magg-X'-first)} for $ t_0 \le \overline t_0$. 
Hence, taking into account \eqref{ex(3-3-magg-X'-first)} and applying Gronwall's lemma to \eqref{ex(3-30-stime-T)}, we deduce the following estimates 
\begin{equation}\label{ex(3-3-Phi-def-1)}
\left| {\Phi \left( {x_d } \right) - \overline \Phi  \left( {x_d } \right)} \right| \le C_{3a} \left[ {\left| {t_0  - \overline t _0 } \right| +\left| {x_0  - \overline x _0 } \right| + \left\| {h - \overline h } \right\|_\infty  } \right], 
\end{equation}
\begin{equation}\label{ex(3-3-X-def-def-1)}
\left| {{\rm X}'\left( {\Phi \left( {x_d } \right);t_0 ,x_0 ,h} \right) - {\rm X'}\left( {\overline \Phi  \left( {x_d } \right);t_0 ,\overline x _0 ,\overline h } \right)} \right| \le C_{3b} \left[ {\left| {t_0  - \overline t _0 } \right| +\left| {x_0  - \overline x _0 } \right| + \left\| {h - \overline h } \right\|_\infty  } \right], 
\end{equation}
for all $x_d \in \left[ {x_{0d}  \vee \overline x _{0d} ,\beta } \right]$, where $C_{3a}$ and $C_{3b}$ are so defined 
\begin{equation}\label{ex(3-1-def-C-3a)}
C_{3a}  = C\left( {1 + \left\| {D_{\left( {x,y} \right)} b} \right\|_{L^1 (0,t_1;h,\overline h )} } \right)\times
\end{equation}
$$
\times \exp \left[ {C\left( {\left\| {b} \right\|_{L^\infty  (0,t_1;h,\overline h )}  + \Lambda \left( {D _x h,D _x \overline h } \right)\left\| {D_{\left( {x,y} \right)} b} \right\|_{L^1 (0,t_1;h,\overline h )} } \right)} \right],
$$
\begin{equation}\label{ex(3-1-def-C-3b)}
C_{3b}  = C\left( {1 + \left\| {b} \right\|_{L^\infty  (0,t_1;h,\overline h )}^2 +\left\| {D_{\left( {x,y} \right)} b} \right\|_{L^1 (0,t_1;h,\overline h )}^2 } \right)\times
\end{equation}
$$
\times \exp \left[ {C\left( {\left\| {b} \right\|_{L^\infty  (0,t_1;h,\overline h )}  + \Lambda \left( {D _x h,D_x \overline h } \right)\left\| {D_{\left( {x,y} \right)} b} \right\|_{L^1 (0,t_1;h,\overline h )} } \right)} \right]
$$
and $C$ depends only on $B_d, \left\| {D _{\left( {t,x'} \right)} b_d } \right\|_{L_{x_d }^1 ( {0,1;L_{(t,x')}^\infty  (S_{t_1}' )} )} $.
 
After having obtained \eqref{ex(3-1-norma-X-second)}, \eqref{ex(3-3-Phi-def-1)} and \eqref{ex(3-3-X-def-def-1)}, we can prove the following properties

A) if $\tau _ -  \left( {t_0 ,x_0 ,h} \right) = 0,\quad {\rm X}_d \left( {\tau _ -  \left( {t_0 ,x_0 ,h} \right);t_0 ,x_0 ,h} \right) < 1$ then there exists $\delta >0$ such that $\tau _ -  \left( { \overline t _0 ,\overline x _0 ;\overline h } \right) = 0$ for $\left| {t_0  - \overline t _0 } \right|,\left| {x_0  - \overline x _0 } \right|,\left\| {h - \overline h } \right\|_{L^\infty  \left( {0,t_1;W^{1,\infty } \left( G  \right)} \right)^p}  < \delta $;

B) if $\tau _ -  \left( {t_0 ,x_0 ,h} \right) > 0$ there exists $\delta >0$ such that ${\rm X}_d \left( {\tau _ -  \left( { \overline t _0 ,\overline x _0 ;\overline h } \right); \overline t _0 ,\overline x _0 ;\overline h } \right) = 1$ for $\left| {t_0  - \overline t _0 } \right| ,\left| {x_0  - \overline x _0 } \right|,\left\| {h - \overline h } \right\|_{L^\infty  \left( {0,t_1;W^{1,\infty } \left( G  \right)} \right)^p}  < \delta $.

We prove, for example, A). If A) is not true, then there exists a sequence $\left\{ {\left({t_{0n} ,x_{0n} ,h_n } \right)}| n \in \mathbb{N} \right\}  \subseteq 
 {S_{t_1}}  \times {L^\infty  \left( {0,t_1;W^{1,\infty } \left( G  \right)} \right)^p}$ such that
\begin{equation}\label{ex(3-2-succ-n-esima-50)}
\left| {t_0  - t_{0n} } \right|,\left| {x_0  - x_{0n} } \right|, \left\| {h - h_n } \right\|_{L^\infty  \left( {0,t_1;W^{1,\infty } \left( G  \right)} \right)^p}  < \frac{1}{n+1}, 
\end{equation}
$$
\quad \tau _ -  \left( {t_{0n},x_{0n} ,h_n } \right) > 0, \quad {\rm X}_d \left( {\tau _ -  \left( {t_{0n},x_{0n} ,h_n } \right);t_{0n},x_{0n} ,h_n } \right) = 1 \quad \forall n \in \mathbb{N}.
$$
We assume, with reference to \eqref{ex(3-3-Phi-def-1)}, $\overline t_{0}=t_{0n}$, $\overline x_{0}=x_{0n}$ and $\overline h = h_n$; therefore, for $n$ big enough, it is possible to take $\beta={\rm X}_d \left( {\tau _ -  \left( {t_{0,} x_0 ,h} \right);t_{0,} x_0 ,h} \right)$. Hence, assuming $x_d=\beta$ and applying \eqref{ex(3-3-Phi-def-1)}, we deduce that  for every $\epsilon >0$ there exists $n_0 \in \mathbb{N}$ such that
\begin{equation}\label{ex(3-2-stima-tempo-n-esima-50)}
\left| {\tau _ -  \left( {t_{0,} x_0 ,h} \right) - {\rm T}\left( {{\rm X}_d \left( {\tau _ -  \left( {t_{0,} x_0 ,h} \right);t_{0,} x_0 ,h} \right);t_{0n} ,x_{0n} ,h_n } \right)} \right| = 
\end{equation}
$$
=  {{\rm T}\left( {{\rm X}_d \left( {0;t_{0,} x_0 ,h} \right);t_{0n} ,x_{0n} ,h_n } \right)}  \le \epsilon \quad \forall n \geq n_0,
$$
therefore we deduce
\begin{equation}\label{ex(3-2-stima-bis-tempo-n-esima-50)}
\tau _ -  \left( {t_{0n} ,x_{0n} ,h_n } \right) \le \varepsilon \quad \forall n \geq n_0,
\end{equation}
from which it follows immediately
\begin{equation}\label{ex(3-2-stima-lim-tempo-n-esima-50)}
\mathop {\lim }\limits_{n \to \infty } \tau _ -  \left( {t_{0n} ,x_{0n} ,h_n } \right) = 0.
\end{equation}
Applying \eqref{ex(3-1-norma-X-second)}, with $s= \tau _ -  \left( {t_{0n} ,x_{0n} ,h_n } \right)$, $\overline t _0  = t_{0n}$, $\overline x _0  = x_{0n}$, $\overline h  = h_n$ we deduce that exists $n'_0 \in \mathbb{N}$, greater than $n_0$, such that
\begin{equation}\label{ex(3-2-stima-lim-X-n-esima-50)}
\left| {{\rm X}_d \left( {\tau _ -  \left( {t_{0n} ,x_{0n} ,h_n } \right);t_0 ,x_0 ,h} \right) - {\rm X}_d \left( {\tau _ -  \left( {t_{0n} ,x_{0n} ,h_n } \right);t_{0n} ,x_{0n} ,h_n } \right)} \right| = 
\end{equation}
$$
=\left| {{\rm X}_d \left( {\tau _ -  \left( {t_{0n} ,x_{0n} ,h_n } \right);t_0 ,x_0 ,h} \right) - 1} \right| < \varepsilon \quad \forall n \geq n'_0. 
$$
Hence, thanks to \eqref{ex(3-2-stima-lim-tempo-n-esima-50)}, \eqref{ex(3-2-stima-lim-X-n-esima-50)} and to the continuity of ${\rm X}_d \left( { \cdot ;t_0 ,x_0 ,h} \right)$, we obtain ${\rm X}_d \left( {\tau _ -  \left( {t_0 ,x_0 ,h} \right);t_0 ,x_0 ,h} \right) = 1$ and this is a contradiction; therefore, we have showed A).

Now, if the hypothesis of A) are verified, then we have
\begin{equation}\label{ex(3-2-stima-tau-60)}
\left| {\tau _ -  \left( {t_0 ,x_0 ,h} \right) - \tau _ -  \left( { \overline t _0 ,\overline x _0 ;\overline h } \right)} \right|  = 0
\end{equation}
for $\left| {t_0  - \overline t _0 } \right|,\left| {x_0  - \overline x _0 } \right|,\left\| {h - \overline h } \right\|_{L^\infty  \left( {0,t_1;W^{1,\infty } \left( G  \right)} \right)^p}  < \delta.$ Therefore, thanks to  \eqref{ex(3-1-norma-X-second-bis)}, we deduce 
\begin{equation}\label{ex(3-2-stima-X-60-new)}
\left| {{\rm X}\left( {\tau _ -  \left( {t_0 ,x_0 ,h} \right);t_0 ,x_0 ,h} \right) - {\rm X}\left( {\tau _ -  \left( { \overline t _0 ,\overline x _0 ;\overline h } \right);\overline t _0 ,\overline x _0 ,\overline h } \right)} \right| \le 
\end{equation}
$$
C_2\Big[ { \left| {t_0  - \overline t _0 } \right|+ \left| {x_0  - \overline x _0 } \right| + \left\| {h - \overline h } \right\|_\infty  }  \Big]
$$
for $\left| {t_0  - \overline t _0 } \right|,\left| {x_0  - \overline x _0 } \right|,\left\| {h - \overline h } \right\|_{L^\infty  \left( {0,t_1;W^{1,\infty } \left( G  \right)} \right)^p}  < \delta.$

Instead, if the hypothesis of B) is verified, then, thanks to \eqref{ex(3-3-Phi-def-1)} and \eqref{ex(3-3-X-def-def-1)}, with $\beta=1$, we obtain
\begin{equation}\label{ex(3-2-stima-bis-tau-60)}
\left| {\tau _ -  \left( {t_0 ,x_0 ,h} \right) - \tau _ -  \left( { \overline t _0 ,\overline x _0 ;\overline h } \right)} \right|  \le C_{3a} \left[ {  \left| {t_0  - \overline t _0 } \right|+\left| {x_0  - \overline x _0 } \right| 
+ \left\| {h - \overline h } \right\|_{L^\infty  (S_{t_1} )} } \right],
\end{equation}
\begin{equation}\label{ex(3-2-stima-bis-X-60-new)}
\left| {{\rm X}'\left( {\tau _ -  \left( {t_0 ,x_0 ,h} \right);t_0 ,x_0 ,h} \right) - {\rm X'}\left( { \tau _ -  \left( { \overline t _0 ,\overline x _0 ;\overline h } \right);\overline t_0 ,\overline x _0 ,\overline h } \right)} \right| \le  
\end{equation}
$$
C_{3b} \left[ {\left| {t_0  - \overline t _0 } \right|+\left| {x_0  - \overline x _0 } \right| + \left\| {h - \overline h } \right\|_\infty  } \right]
$$
for $\left| {t_0  - \overline t _0 } \right|,\left| {x_0  - \overline x _0 } \right|,\left\| {h - \overline h } \right\|_{L^\infty  \left( {0,t_1;W^{1,\infty } \left( G  \right)} \right)^p}  < \delta.$

Of course it is possible that A) and B) are not verified; in this case, we deduce that $\tau _ -  \left( {t_0 ,x_0 ,h} \right) = 0$ and ${\rm X}_d \left( {\tau _ -  \left( {t_0 ,x_0 ,h} \right),t_0 ,x_0 ,h} \right) = 1$, therefore we can apply \eqref{ex(3-2-stima-tau-60)}-\eqref{ex(3-2-stima-X-60-new)}  or \eqref{ex(3-2-stima-bis-tau-60)}-\eqref{ex(3-2-stima-bis-X-60-new)}. 

Hence, we have showed \eqref{ex(3-1-norma-tau-first)}.

Now, after having fixed $s \in \left( {\tau _ -  \left( {t_0 ,x_0 ,h} \right),t_0 } \right)$ and  using the regularity of $\tau _ -  \left(  \cdot  \right)
$, it is possible to find $\delta>0$ such that
\begin{equation}\label{ex(3-2-magg-X-finale-1)}
\left| {{\rm X}\left( {s;t_0 ,x_0 ,h} \right) - {\rm X}\left( {\overline s ;\overline t _0 ,\overline x _0 ,\overline h } \right)} \right| \le \left| {x_0  - \overline x _0 } \right| + \int\limits_{I(s , \overline s ) \cup I(t_0  , \overline t _0 )} {\left\| {b\left( {r, \cdot } \right)} \right\|_{L^\infty  ({G};h, \overline h )} } dr + \Delta
\end{equation}
if $\left| {t_0  - \overline t _0 } \right|,\left| {x_0  - \overline x _0 } \right|,\left\| {h - \overline h } \right\|_{W^{1,\infty } (S_{t_1} )^p} ,\left| {s - \overline s } \right| \le \delta   
$ and $\overline s \in \left( {\tau _ -  \left( {\overline t_0 ,\overline x_0 ,\overline h} \right),\overline t_0 } \right)$, where
\begin{equation}\label{ex(3-2-def-delta-first)}
\Delta = \int\limits_{s \vee \overline s }^{t_0  \wedge \overline t _0 } {\left| {b\left( {r,{\rm X}\left( {r;t_0 ,x_0 ,h} \right),h\left( {r,{\rm X}\left( {r;t_0 ,x_0 ,h} \right)} \right)} \right) - } \right.} 
\end{equation}
$$
\left. { - b\left( {r,{\rm X}\left( {r;\overline t _0 ,\overline x _0 ,\overline h } \right), \overline h\left( {r,{\rm X}\left( {r;\overline t _0 ,\overline x _0 ,\overline h } \right)} \right)} \right)} \right|dr.
$$
After simple calculations, we have 
\begin{equation}\label{ex(3-2-stima-delta-first)}
\Delta \le \int\limits_{s \vee \overline s }^{t_0  \wedge \overline t _0 } {\left\| {D_{\left( {x,y} \right)} b} \right\|_{L^\infty  ({G};h,\overline h )} } \left[ {\Lambda \left( {D_x h,D_x \overline h } \right)
   \times } \right.
\end{equation}
$$
\left. { \times \left| {{\rm X}\left( {r;t_0 ,x_0 ,h} \right) - {\rm X}\left( {r;\overline t _0 ,\overline x _0 ,\overline h } \right)} \right| + \left\| {h - \overline h } \right\|_\infty  } \right]dr.
$$
Therefore, we can obtain \eqref{ex(4-1-terza-stima-X)} by \eqref{ex(3-2-magg-X-finale-1)}, \eqref{ex(3-2-stima-delta-first)} and \eqref{ex(3-1-norma-X-second)}.

\  $ \square$ 
\medskip
\medskip


\medskip 
\medskip 
Now we are ready to study the continuous dependence of flow $\rm X$ and initial time $\tau_-$ on vector field
$b$ and parametric vector function $h$.

{\bf Lemma 5.3.}   \ {\it If $b^{\left( k \right)} \in L_t^1 \big( {0,t_1;W_{(x,y_{loc} )}^{1,\infty } \left( {{G}  \times \mathbb{R}^p} \right)} \big)^d$, $B_d^{\left( k \right)} >0$ such that $b_d^{(k)} (t,x) \le  - B_d^{(k)} 
$ $\forall x \in G$ a.e. $t \in (0,t_1)$, $(t_0,x_0,h^{\left( k \right)}) \in$ $S_{t_1} \times$ $L^\infty  \left( {0,T;W^{1,\infty } \left( G  \right)} \right)^p
$;  ${\rm X}^{\left( k \right)} \left( { \cdot ;t_0 ,x_0 ,h^{\left( k \right)}} \right):[ {\tau^{\left( k \right)} _ -  (t_0 ,x_0 ,h^{\left( k \right)}),t_0 } ] \to \mathbb{R}^d$  is the flux associated with $b^{\left( k \right)}$ where $k=1,2$,  then we have 
\begin{equation}\label{ex(4-2-n-stima-diff-X)}
\left| {{\rm X}^{\left( 2 \right)}\left( {s;t_0 ,x_0 ,h^{\left( 2 \right)}} \right) -  {\rm X}^{\left( 1 \right)} \left( {s;t_0 ,x_0 , h^{\left( 1 \right)} } \right)} \right| \le
\end{equation}
$$
\le C_1 \Big( {\left\| {h^{\left( 2 \right)} - h^{\left( 1 \right)} } \right\|_\infty   + \left\| {b^{\left( 2 \right)} - b^{\left( 1 \right)} } \right\|_{L^1 (0,t_1;h^{\left( 1 \right)},h^{\left( 2 \right)})} } \Big),
$$
for every $s \in [ \tau _ - ^{(1)} 
  (t_0 ,x_0 ,h^{\left( 1 \right)}) \vee  \tau _ - ^{(2)} 
  \left( {t_0 ,x_0 , h^{\left( 2 \right)} } \right),t_0 ]
$; furthermore, the constant $C_1$ is determined by assuming $b=b^{(1)}$, $h=h^{(1)}$ and $\overline h=h^{(2)}$.

Finally, if $b^{\left( k \right)} \in  L_{x_d}^{-} \left( {S_{t_1} \times \mathbb{R}^p } \right)
$ (see \eqref{ex(4-1-spazi-4-12)})  with $k=1,2$, then for every $M>0$ there exists $\delta >0$ depending on $h^{\left( 1 \right)},b^{\left( 1 \right)},M$, such that, under the conditions $\left\| {h^{\left( 2 \right)} - h^{\left( 1 \right)} } \right\|_{\infty} < \delta$, ${\left\| {b^{\left( 2 \right)} - b^{\left( 1 \right)} } \right\|_{L^1 (0,t_1;h^{\left( 1 \right)},h^{\left( 2 \right)})} }$
$< \delta$ and $\left\| {b^{(k)} } \right\|_{L^\infty  (0,t_1 ;h^{(k)} )}  \le M
$ with $k=1,2$, we have
\begin{equation}\label{ex(4-2-n-stima-tau)}
\left| {\tau _ - ^{(2)}  (t_0 ,x_0 ,h^{\left( 2 \right)}) - \tau _ - ^{(1)}  (t_0 ,x_0 , h^{\left( 1 \right)})} \right| + 
\end{equation}
$$
+\left| {{\rm X}^{\left( 2 \right)}(\tau _ - ^{(2)}   (t_0 ,x_0 ,h^{\left( 2 \right)});t_0 ,x_0 ,h^{\left( 2 \right)}) - {\rm X}^{\left( 1 \right)}( \tau _ - ^{(1)}  (t_0 ,x_0 , h^{\left( 1 \right)} );t_0 ,x_0 ,h^{\left( 1 \right)})} \right|
\le 
$$
$$
\le C_5 \Big( {\left\| {h^{\left( 2 \right)} -  h^{\left( 1 \right)} } \right\|_\infty   + \left\| {b^{\left( 2 \right)} - b^{\left( 1 \right)} } \right\|_{L^1 (0,t_1;h^{\left( 1 \right)}, h^{\left( 2 \right)})} } \Big),
$$
where 
\begin{equation}\label{ex(4-1-def-c-8)}
C_5  = C \Big( {1 + \mathop  \bigvee \limits_{k = 1,2} \left\| {b^{(k)} } \right\|_{L^\infty  (0,t_1 ;h^{(k)}  )} } \Big) \Big( 1 + \left\| {D_{(x,y)} b^{(1)} } \right\|_{L^1 (0,t_1 ;h^{(1)} ,h^{(2)} )} 
 \Big) \times 
\end{equation}
$$
\times \exp \Big[ C \Big( { \mathop  \bigvee \limits_{k = 1,2} \left\| {b^{(k)} } \right\|_{L^\infty  (0,t_1 ;h^{(k)} )} } + \Lambda \left( {D_x h^{(1)} ,D_x h^{(2)} } \right)\left\| {D_{(x,y)} b^{(1)} } \right\|_{L^1 (0,t_1 ;h^{(1)} ,h^{(2)} )} 
   \Big)\Big]$$
with the constant $C$  depending only on  $B_d^{\left( k \right)}$, $\left\|{  b_d^{\left( k \right)} } \right\|_{L_{x_d }^1 (0,1;L_{(t,x')}^\infty  (S_{t_1}'  ))} $ $($with $k=1,2)$,  $\left\| {D _{\left( {t,x'} \right)} b_d^{\left( 1 \right)} } \right\|_{L_{x_d }^1 (0,1;L_{(t,x')}^\infty  (S_{t_1}'  ))} $.

}

\medskip

{\sc {\bf Proof}}. \

Let us define
\begin{equation}\label{ex(4-2-k-equaz-integrale)}
 { \rm X}^{\left( k \right)}\left( {s;t_0 ,x_0 ; h^{\left( k \right)}} \right) = x_0  - \int\limits_s^{t_0 } { b^{\left( k \right)}\left( {r, {\rm X}^{\left( k \right)}\left( {r;t_0 ,x_0 ; h^{\left( k \right)}} \right); h^{\left( k \right)}\left( {r, {\rm X}^{\left( k \right)}\left( {r;t_0 ,x_0 ; h^{\left( k \right)}} \right)} \right)
} \right)dr}, 
\end{equation}
where $s \in [ \tau _ - ^{\left( k \right)} \left( {t_0 ,x_0 ,h^{\left( k \right)} )},t_0 \right]$ and $k=1,2$. We obtain \eqref{ex(4-2-n-stima-diff-X)} by the difference between \eqref{ex(4-2-k-equaz-integrale)} with $k=2$ and \eqref{ex(4-2-k-equaz-integrale)} with $k=1$.

Now, using similar arguments to those employed to obtain the equation \eqref{ex(3-30-def-eq-int-T)}, we can transform the equation \eqref{ex(4-2-k-equaz-integrale)} in the following vectorial integral equation
\begin{equation}\label{ex(4-2-n-def-eq-int-T)}
\Phi^{\left( k \right)}\left( {x_d } \right) = t_0  + \int\limits_{x_{0d} }^{x_d } {\frac{{dw}}{{ b_d^{\left( k \right)} \left( { \Phi^{\left( k \right)} \left( w \right), { \rm X}'^{\left( k \right)}\left( { \Phi^{\left( k \right)} \left( w \right);Z^{\left( k \right)}_0} \right),w} \right)
 }}}, 
\end{equation}
\begin{equation}\label{ex(4-2-n-def-eq-int-X')}
 {\rm X}'^{\left( k \right)}( \Phi^{\left( k \right)} \left( {x_d } \right) ;Z_0^{(k)}) = x_0 ' + 
\end{equation}
$$
+\int\limits_{x_{0d} }^{x_d } {\frac{{ b'^{\left( k \right)}( \Phi^{\left( k \right)} \left( {w } \right) ,X^{'(k)} (\Phi ^{(k)} (w);Z_0^{(k)} )
,w, h^{(k)} (\Phi ^{(k)} (w),X^{'(k)} (\Phi ^{(k)} (w);Z_0^{(k)} ),w)
) }}{{ b_d^{\left( k \right)} ( \Phi^{\left( {k} \right)} \left( {w } \right) , X^{'(k)} (\Phi ^{(k)} (w);Z_0^{(k)} )
,w) }}dw}, 
$$
for every $x_d  \in [  x_{0d} , {\rm X}_d^{\left( k \right)} ( \tau _ - ^{(k)}  (t_0 ,x_0 , h^{\left( k \right)});t_0 ,x_0 , h^{\left( k \right)})]$, $k=1,2$, where we have defined 
\begin{equation}\label{ex(4-2-n-def-T)}
Z_0^{(k)}=(t_0 ,x_0 , h^{\left( k \right)}) ,\quad \Phi^{\left( k \right)}  \left(  \cdot  \right) = {\rm T}^{\left( k \right)} (  \cdot ;Z_0^{(k)} ),
\end{equation}
and ${\rm T}^{\left( k \right)} (  \cdot ;Z_0^{(k)} )$ is the inverse function of ${\rm X}^{\left( k \right)} (  \cdot ;Z_0^{(k)} )$.

Studying the difference between \eqref{ex(4-2-n-def-eq-int-T)} with $k=2$ and \eqref{ex(4-2-n-def-eq-int-T)} with $k=1$, we can deduce the following estimate
\begin{equation}\label{ex(4-2-n-diff-T-1)}
\left| {\Phi ^{(2)} \left( {x_d } \right) - \Phi ^{(1)} \left( {x_d } \right)} \right| \le B_d^{(2) - 1} B_d^{(1) - 1} \Big[\left\| {b_d^{(2)}  - b_d^{(1)} } \right\|_{L_{x_d }^1 (0,1;h^{(1)} ,h^{(2)} )}  +
\end{equation}
$$
+\int\limits_{x_{0d} }^{x_d } {\left\| {D_{(t,x')} b_d^{(1)} ( \cdot ,w)} \right\|_{L^\infty  (S'_{t_1 } )} \big(\left| {\Phi ^{(2)} \left( w \right) - \Phi ^{(1)} \left( w \right)} \right|}  + 
$$
$$
+\big| {{\rm X^{'(2)}}(\Phi ^{(2)} \left( w \right);Z_0^{(2)} ) - {\rm X^{'(1)}}(\Phi ^{(1)} \left( w \right);Z_0^{(1)} )} \big|\big) dw\Big],
$$
where $x_d \in  [  x_{0d} , \mathop  \bigwedge \limits_{k = 1,2}{\rm X}_d^{\left( k \right)} ( \tau _ - ^{(k)}  (t_0 ,x_0 , h^{\left( k \right)});t_0 ,x_0 , h^{\left( k \right)})].$ 

Now, thanks to \eqref{ex(3-1-norma-X-first-bis)} and \eqref{ex(4-2-n-stima-diff-X)}, we have
\begin{equation}\label{ex(4-2-n-diff-X'-1)}
\big| {{\rm X^{'(2)}}(\Phi ^{(2)} \left( w \right);Z_0^{(2)} ) - {\rm X^{'(1)}}(\Phi ^{(1)} \left( w \right);Z_0^{(1)} )} \big| \le \mathop  \bigvee \limits_{k = 1,2} \left\| {b^{(k)} } \right\|_{L^\infty  (0,t_1 ;h^{(k)} )} \left| {\Phi ^{(2)} \left( w \right) - \Phi ^{(1)} \left( w \right)} \right| + 
\end{equation}
$$
+C_1 \left( {\left\| {h^{\left( 2 \right)} -  h^{\left( 1 \right)} } \right\|_\infty   + \left\| {b^{\left( 2 \right)} -  b^{\left( 1 \right)} } \right\|_{L^1 (0,t_1;h^{\left( 1 \right)}, h^{\left( 2 \right)})} } \right), \quad \forall w \in [x_{0d},x_d].
$$
Hence, applying \eqref{ex(4-2-n-diff-X'-1)} in \eqref{ex(4-2-n-diff-T-1)} and Gronwall's lemma, we deduce 
\begin{equation}\label{ex(4-2-n-diff-T-2)}
\left| {{\rm T}^{\left( 2 \right)}\left( {x_d ;t_0 ,x_0 ,h^{\left( 2 \right)}} \right) - {\rm T}^{\left( 1 \right)}\left( { x_d ;t _0 , x _0 ; h^{\left( 1 \right)} } \right) } \right| \le 
\end{equation}
$$
\le C'_{3a }\left( {\left\| {h^{\left( 2 \right)} -  h^{\left( 1 \right)} } \right\|_\infty   + \left\| {b^{\left( 2 \right)} -  b^{\left( 1 \right)} } \right\|_{L^1 (0,t_1;h^{\left( 1 \right)}, h^{\left( 2 \right)})} } \right),
$$
\begin{equation}\label{ex(4-2-n-diff-X'-2)}
\left| {{\rm X}'^{\left( 2 \right)}\left( {{\rm T}^{\left( 2 \right)}\left( { x_d ;t_0 ,x_0 ,h^{\left( 2 \right)}} \right);t_0 ,x_0 ,h^{\left( 2 \right)}} \right) - {\rm X}'^{\left( 1 \right)}\left( {{\rm T}^{\left( 1 \right)}\left( { x_d ; t _0 ,x _0 ;h^{\left( 1 \right)} } \right);t _0 ,x _0 ;h^{\left( 1 \right)} } \right)} \right| \le
\end{equation}
$$
\le C'_{3b }\left( {\left\| {h^{\left( 2 \right)} -  h^{\left( 1 \right)} } \right\|_\infty   + \left\| {b^{\left( 2 \right)} -  b^{\left( 1 \right)} } \right\|_{L^1 (0,t_1;h^{\left( 1 \right)}, h^{\left( 2 \right)})} } \right),
$$
where
\begin{equation}\label{ex(4-2-def-c'3a)}
C'_{3a }  = C \Big( 1 + \left\| {D_{(x,y)} b^{(1)} } \right\|_{L^1 (0,t_1 ;h^{(1)} ,h^{(2)} )} 
 \Big) \times 
\end{equation}
$$
\times \exp \Big[ C \Big( { \mathop  \bigvee \limits_{k = 1,2} \left\| {b^{(k)} } \right\|_{L^\infty  (0,t_1 ;h^{(k)} )} } + \Lambda \left( {D _x h^{(1)} ,D_x h^{(2)} } \right)\left\| {D_{(x,y)} b^{(1)} } \right\|_{L^1 (0,t_1 ;h^{(1)} ,h^{(2)} )} 
   \Big)\Big],$$
\begin{equation}\label{ex(4-2-def-c'3b)}
C'_{3b }=\Big( {1 + \mathop  \bigvee \limits_{k = 1,2} \left\| {b^{(k)} } \right\|_{L^\infty  (0,t_1 ;h^{(k)}  )} } \Big)C'_{3a }.
\end{equation}

After having obtained \eqref{ex(4-2-n-stima-diff-X)}, \eqref{ex(4-2-n-diff-T-2)},  \eqref{ex(4-2-n-diff-X'-2)} and assuming $M>0$, we can prove the following properties

$A' )$ if $\tau _ - ^{\left( 1 \right)} (t_0 ,x_0 ,h^{\left( 1 \right)} ) = 0,\quad {\rm X}_d ^{\left( 1 \right)}(\tau _ - ^{\left( 1 \right)} (t_0 ,x_0 ,h^{\left( 1 \right)} );t_0 ,x_0 ,h^{\left( 1 \right)}) < 1$ then there exists $\delta >0$ such that $\tau _ - ^{\left( 2 \right)} (t_0 ,x_0 ,h^{\left( 2 \right)} )= 0$ for $\left\| {h^{\left( 1 \right)} -  h^{\left( 2 \right)} } \right\|_{\infty}$, ${\left\| {b^{\left( 1 \right)} - b^{\left( 2 \right)} } \right\|_{L^1 (0,t_1;h^{\left( 1 \right)}, h^{\left( 2 \right)})} }
  < \delta $ and $\left\| {b^{(k)} } \right\|_{L^\infty  (0,t_1 ;h^{(k)} )}  \le M
$ with $k=1,2$;

$B' )$ if $\tau _ - ^{\left( 1 \right)} (t_0 ,x_0 ,h^{\left( 1 \right)} ) > 0$ there exists $\delta$ $>0$ such that ${\rm X}_d ^{\left( 2 \right)}$$(\tau _ - ^{\left( 2 \right)} (t_0 ,x_0 ,h^{\left( 2 \right)} );t_0 ,x_0 ,h^{\left( 2 \right)})$$= 1$ for $\left\| {h^{\left( 1 \right)} -  h^{\left( 2 \right)} } \right\|_{\infty}$, ${\left\| {b^{\left( 1 \right)} - b^{\left( 2 \right)} } \right\|_{L^1 (0,t_1;h^{\left( 1 \right)}, h^{\left( 2 \right)})} }$
  $< \delta $ and $\left\| {b^{(k)} } \right\|_{L^\infty  (0,t_1 ;h^{(k)} )}  \le M
$ with $k=1,2$.

Using \eqref{ex(4-2-n-stima-diff-X)}, \eqref{ex(4-2-n-diff-T-2)} and \eqref{ex(4-2-n-diff-X'-2)}, the proofs of $ A' )$,$ B' )$ may be carried out as those used to show A) and B) (see the proof of  Lemma 5.2).  At this point, to show \eqref{ex(4-2-n-stima-tau)}, we use the same argument employed in the proof of \eqref{ex(3-1-norma-tau-first)}.  

\  $ \square$ 
\medskip
\medskip
\medskip
\medskip
\section{- The well-posedness of the IBVP for the linear transport equation.}

Before studying the initial and boundary value problem  \eqref{ex(4-1-lin-transp-pde-mixed)}-\eqref{ex(4-1-condiz-boundary-value)}, we give the definition of generalized solution for this problem following the one given in \cite{[M1]}.

{\bf Definition 6.0.}   \ {
\it Let us assume $h \in L^\infty  \left( {0,t_1 ;W^{1,\infty } \left( G  \right)} \right)^p$, \eqref{ex(3-1-condiz-b-a)}-\eqref{ex(3-1-condiz-bd)} for the vector field $b$ and  
\begin{equation}\label{ex(4-2-2-reg-z-*-)}
z^*  \in C^0 \left( {\Gamma _ -  } \right) \cap L^\infty  \left( {\Gamma _ -  } \right),
\end{equation}
\begin{equation}\label{ex(4-2-2-reg-c-f-0)}
c,a \in L^\infty_t \big( 
{0,t_1;L_{(x,y_{loc} )}^{\infty } \left( {G  \times \mathbb{R}^p} \right)} \big)  \cap  L^1_t \big( 
{0,t_1;W_{(x,y_{loc} )}^{1,\infty } \left( {G  \times \mathbb{R}^p} \right)} \big).
\end{equation}

Hence, we say that the function 
$z(\cdot,h) \in L^\infty  \left( {0,t_1 ;W^{1,\infty } \left( G  \right)} \right)$ is a generalized solution for the IBVP \eqref{ex(4-1-lin-transp-pde-mixed)}-\eqref{ex(4-1-condiz-boundary-value)}, if and only if
\begin{equation}\label{ex(4-2-2-int-eq)}
z\left( {t,x;h} \right) = z^* \left( {\tau _ -  \left( {t,x,h} \right),{\rm X}\left( {\tau _ -  \left( {t,x,h} \right);t,x,h} \right)} \right) +  
\end{equation}
$$
+\int\limits_{\tau _ -  \left( {t,x,h} \right)}^t {\left[ -{c\left( {s,{\rm X}\left( {s;t,x,h} \right),h\left( {s,{\rm X}\left( {s;t,x,h} \right)} \right)} \right)} \right.} z\left( {s,{\rm X}\left( {s;t,x,h} \right);h} \right) + 
$$
$$
 + \left. {a\left( {s,{\rm X}\left( {s;t,x,h} \right),h\left( {s,{\rm X}\left( {s;t,x,h} \right)} \right)} \right)} \right]ds
$$
for every $(t,x) \in S_{t_1}$. 
}

Now, we are ready to give some results on the existence, uniqueness, regularity and continuous dependence on data of the generalized solution for the problem \eqref{ex(4-1-lin-transp-pde-mixed)}-\eqref{ex(4-1-condiz-boundary-value)}.

{\bf Lemma 6.1.}   \ {\it Let us assume  $b \in L_{x_d }^ -  \left( {S_{t_1 } } \right)
$,  $z^* \in Lip^{unif}_{loc}\left( {\Gamma _ -  } \right)$  and \eqref{ex(4-2-2-reg-c-f-0)} for $c,$ $a$. Furthermore, we denote by $L^*$ a positive constant such that 
\begin{equation}\label{ex(4-2-2-lipschitz-z-*)}
\left| {z^* \left( {t,x} \right) - z^* \left( {\overline t ,\overline x } \right)} \right| \le L^* \left\{ {\left| {t - \overline t } \right| + \left| {x - \overline x } \right|} \right\},
\end{equation}
where $\left( {t,x} \right),\left( {\overline t ,\overline x } \right) \in \Gamma _ - $ and there exists a positive number $\sigma>0$ such that $\left| {t - \overline t} \right| + \left| {x - \overline x} \right| \le \sigma. 
$

Then, for every $h \in L^\infty  \left( {0,t_1 ;W^{1,\infty } \left( G  \right)} \right)^p$ there exists one and only one generalized solution $z(\cdot,h) \in L^\infty  \left( {0,t_1 ;W^{1,\infty } \left( G  \right)} \right)$ for IBVP \eqref{ex(4-1-lin-transp-pde-mixed)}-\eqref{ex(4-1-condiz-boundary-value)}.

Furthermore, this solution satisfies the following inequality 
\begin{equation}\label{ex(4-2-2-stima-z)}
\left\| {z\left( { \cdot ,h} \right)} \right\|_{L^\infty  ( S_{t_1} )}  \le \left[ {\left\| {z^* } \right\|_{L^\infty  \left( {\Gamma _ -  } \right)}  + \left\| a \right\|_{L^1 \left( {0,t_1;h} \right)} } \right]\exp  {\left\| c \right\|_{L^1 \left( {0,t_1;h} \right)} } .
\end{equation}

}

\medskip

{\sc {\bf Proof}}. \

First of all, we show the uniqueness of the generalized solution. Hence, we suppose that $z(\cdot,h)$ and  $w(\cdot,h)$ are two generalized solutions for the IBVP \eqref{ex(4-1-lin-transp-pde-mixed)}-\eqref{ex(4-1-condiz-boundary-value)}. Therefore, we immediately deduce
\begin{equation}\label{ex(4-2-2-stima-unicita)}
\left| {z\left( {t,x;h} \right) - w\left( {t,x;h} \right)} \right| \le 
\end{equation}
$$
 \le \left\| c \right\|_{L^\infty  \left( {0,t_1;h } \right)} \int\limits_{\tau _ -  \left( {t,x,h} \right)}^t {\left| {z\left( {s,{\rm X}\left( {\tau _ -  \left( {t,x,h} \right);t,x,h} \right);h} \right) - w\left( {s,{\rm X}\left( {\tau _ -  \left( {t,x,h} \right);t,x,h} \right)} \right)} \right|} ds.
$$
After simple calculations and applying Gronwall's lemma we obtain $z=w$.

Using Gronwall's lemma, the estimate \eqref{ex(4-2-2-stima-z)} follows immediately.

Now, we prove the existence of the generalized solution; for this purpose, we can use the method of successive approximations to prove the existence of generalized solution. Assuming $z_1(\cdot,h)  \in L^\infty  \left( {0,t_1 ;W^{1,\infty } \left( G  \right)} \right)$, we consider the following recurrence schema 
\begin{equation}\label{ex(4-2-2-int-eq-n)}
z_{n+1}\left( {t,x;h} \right) = z^* \left( {\tau _ -  \left( {t,x,h} \right),{\rm X}\left( {\tau _ -  \left( {t,x,h} \right);t,x,h} \right)} \right) +  
\end{equation}
$$
+\int\limits_{\tau _ -  \left( {t,x,h} \right)}^t {\left[ -{c\left( {s,{\rm X}\left( {s;t,x,h} \right),h\left( {s,{\rm X}\left( {s;t,x,h} \right)} \right)} \right)} \right.} z_n \left( {s,{\rm X}\left( {s;t,x,h} \right);h} \right) + 
$$
$$
 + \left. {a\left( {s,{\rm X}\left( {s;t,x,h} \right),h\left( {s,{\rm X}\left( {s;t,x,h} \right)} \right)} \right)} \right]ds, \quad \forall n > 1.
$$ 
We must verify that $z_{n+1}(\cdot;h)$ is well defined and belongs to $L^\infty  \left( {0,t_1 ;W^{1,\infty } \left( G  \right)} \right)$. To show this, it is sufficient to see that if we assume the following inductive hypothesis: $z_{2}(\cdot;h),...,z_{n}(\cdot;h) \in L^\infty  \left( {0,t_1 ;W^{1,\infty } \left( G  \right)} \right)$, then $z_{n+1}(\cdot;h) \in L^\infty  \left( {0,t_1 ;W^{1,\infty } \left( G  \right)} \right)$.

It is not so hard to obtain the following estimates:
\begin{equation}\label{ex(4-2-2-stima-zn+1-zn)}
\left\| {z_{n + 1} \left( { \cdot ,h} \right) - z_n \left( { \cdot ,h} \right)} \right\|_{L^\infty  (S_{t_1} )}  \le \left\| c \right\|_{L^\infty  (0,1;h )}^{n - 1} \left\| {z_2 \left( { \cdot ,h} \right) - z_1 \left( { \cdot ,h} \right)} \right\|_{L^\infty  (S_{t_1} )} \frac{{t^{n - 1}_1 }}{{\left( {n - 1} \right)!}},
\end{equation}
\begin{equation}\label{ex(4-2-2-stima-zn+1)}
\left\| {z_{n + 1} ( \cdot ;h)} \right\|_{L^\infty  (s_{t_1 } )}  \le \left\| {z^* } \right\|_{L^\infty  (\Gamma _ -  )} \sum\limits_{k = 0}^{n - 1} {\left\| c \right\|_{L^\infty  (0,t_1 ;h)}^k \frac{{t_1^k }}{{k!}} + } \left\| a \right\|_{L^\infty  (0,t_1 ;h)} \sum\limits_{k = 0}^{n - 1} {\left\| c \right\|_{L^\infty  (0,t_1 ;h)}^k \frac{{t_1^{k + 1} }}{{(k + 1)!}} + } 
\end{equation}
$$
+ \left\| c \right\|_{L^\infty  (0,t_1 ;h)}^n \left\| {z_1 ( \cdot ;h)} \right\|_{L^\infty  (s_{t_1 } )} \frac{{t_1 ^n }}{{n!}}.
$$
After these preliminary estimates, we are ready to study the regularity of ${z_{n + 1} ( \cdot ;h)}$. Assuming $(t,x),(t,\overline x) \in S_{t_1 }$ and remembering \eqref{ex(4-2-2-int-eq-n)}, we deduce
\begin{equation}\label{ex(4-2-2-stima-diff-zn+1)}
\left| {z_{n + 1} (t,x;h) - z_{n + 1} (t,\overline x ;h)} \right| \le \Delta ^*  + \Delta _a  + \Delta _c .
\end{equation}
$\Delta ^*  , \Delta _a  , \Delta _c$ are so defined 
\begin{equation}\label{ex(4-2-2-stima-delta-*)}
\Delta ^* =\left| {z^* (\tau _ -  (t,x,h),{\rm X}(\tau _ -  (t,x,h);t,x,h)) - z^* (\tau _ -  (t,\overline x ,h),{\rm X}(\tau _ -  (t,\overline x ,h);t,\overline x ,h))} \right| \le
\end{equation}
$$
\le L^* \big( \left| {\tau _ -  (t,x,h) - \tau _ -  (t,\overline x ,h)} \right| + \left| {{\rm X}(\tau _ -  (t,x,h);t,x,h) - {\rm X}(\tau _ -  (t,\overline x ,h);t,\overline x ,h)} \right|
    \big) \le
$$
$$
\le L^* C_3 \left| {x-\overline x} \right|,
$$
where the last inequality can be deduced applying \eqref{ex(3-1-norma-tau-first)} and assuming $x_0=x$, $\overline x_0 = \overline x$, $\overline h =h$ and $\left| {x-\overline x} \right| \le \delta(t,x,h)$;
\begin{equation}\label{ex(4-2-2-stima-delta-a)}
\Delta_a  =\left| {\int\limits_{\tau _ -  (t,x,h)}^t {a\left( {s,{\rm X}(s;t,x,h),h(s;{\rm X}(s;t,x,h))} \right)} ds} \right. - 
\end{equation}
$$
\left. { - \int\limits_{\tau _ -  (t,\overline x ,h)}^t {a\left( {s,{\rm X}(s;t,\overline x ,h),h(s;{\rm X}(s;t,\overline x ,h))} \right)} ds} \right| \le \int\limits_{I(\tau _ -  (t,x,h),\tau _ -  (t,\overline x ,h))} {\left\| {a(s, \cdot )} \right\|_{L^\infty  (G;h)} ds}  + 
$$
$$
+\int\limits_{\tau _ -  (t,x,h) \vee \tau _ -  (t,\overline x ,h)}^t {\left\| {\nabla _{(x,y)} a(s, \cdot )} \right\|_{L^\infty  (G;h)} \big[\left| {{\rm X}(s;t,x,h) - {\rm X}(s;t,\overline x ,h)} \right|}  + 
$$
$$
 + \left| {h(s;{\rm X}(s;t,x,h)) - h(s;{\rm X}(s;t,\overline x ,h))} \right|\big]ds \le
$$
$$
\le \left( {\Lambda \left( {D_x h,D_x h} \right)C_2 \left\| {\nabla _{(x,y)} a} \right\|_{L^1 (0,t_1 ;h)}  + C_3\left\| a \right\|_{L^\infty  (0,t_1 ;h)} } \right)\left| {x - \overline x } \right|,
$$
where the last inequality has been obtained applying \eqref{ex(3-1-norma-X-second-bis)}, \eqref{ex(3-1-norma-tau-first)} and assuming $\left| {x-\overline x} \right| \le \delta(t,x,h)$;
\begin{equation}\label{ex(4-2-2-stima-delta-c)}
\Delta_c  = \left| {\int\limits_{\tau _ -  (t,x,h)}^t {c\left( {s,{\rm X}(s;t,x,h),h(s;{\rm X}(s;t,x,h))} \right)} z_n \left( {s,{\rm X}(s;t,x,h);h)} \right)ds} \right. - 
\end{equation}
$$
\left. { - \int\limits_{\tau _ -  (t,\overline x ,h)}^t {c\left( {s,{\rm X}(s;t,\overline x ,h),h(s;{\rm X}(s;t,\overline x ,h))} \right)} z_n \left( {s,{\rm X}(s;t,\overline x ,h);h)} \right)ds} \right| \le
$$
$$
\le \left\| {z_n \left( { \cdot ;h} \right)} \right\|_{L^\infty  (S_{t_1 } )} \left( {\Lambda \left( {D_x h,D_x h} \right)C_2 \left\| {\nabla _{(x,y)} c} \right\|_{L^1 (0,t_1 ;h)}  + C_3 \left\| c \right\|_{L^\infty  (0,t_1 ;h)} } \right)\left| {x - \overline x } \right| + 
$$
$$
 + \left\| c \right\|_{L^\infty  (0,t_1 ;h)} \int\limits_0^t {\left| {z_n (s,x;h) - z_n (s,\overline x ;h)} \right|}ds,\quad \quad \mbox{where} \left| {x-\overline x} \right| \le \delta(t,x,h).
$$
Therefore, taking into account \eqref{ex(4-2-2-stima-zn+1)}-\eqref{ex(4-2-2-stima-delta-c)}, we deduce there exists a constant $C>0$ such that
\begin{equation}\label{ex(4-2-2-stima-diff-zn+1-bis)}
\left| {z_{n + 1} (t,x;h) - z_{n + 1} (t,\overline x ;h)} \right| \le C\left| {x - \overline x } \right| +\left\| c \right\|_{L^\infty  (0,t_1 ;h)} \int\limits_0^t {\left| {z_n (s,x;h) - z_n (s,\overline x ;h)} \right|}ds,
\end{equation}
where $\left| {x-\overline x} \right| \le \delta(t,x,h)$. Hence $z_{n + 1}(\cdot,h) \in L^\infty  \left( {0,t_1 ;W^{1,\infty } \left( G  \right)} \right)$ and, moreover, thanks to \eqref{ex(4-2-2-int-eq-n)}, is defined on all $S_{t_1}$  .

Now, using \eqref{ex(4-2-2-stima-zn+1-zn)}, we observe that the telescopic series $z_1(\cdot,h) +\sum\limits_{n = 1}^{ + \infty } {\left( {z_{n + 1}(\cdot,h)  - z_n(\cdot,h) } \right)}$ is uniform convergent on $S_{t_1}$ to a function that we denote by $z(\cdot,h)$. Finally, it is not so hard to verify, taking the limit for $n  \to \infty$ in \eqref{ex(4-2-2-stima-diff-zn+1-bis)} and \eqref{ex(4-2-2-int-eq-n)}, that $z(\cdot,h) \in L^\infty  \left( {0,t_1 ;W^{1,\infty } \left( G  \right)} \right)$ and it is the generalized solution for the IBVP \eqref{ex(4-1-lin-transp-pde-mixed)}-\eqref{ex(4-1-condiz-boundary-value)}.

\  $ \square$ 
\medskip
\medskip

{\bf Lemma 6.2.}   \ {\it Let us assume the hyphoteses of Lemma 6.1.  

For every $\left( {t,x,h} \right) \in S_{t_1}  \times L^\infty  \left( {0,t_1 ;W^{1,\infty } \left( G  \right)} \right)^p$ there exists $\delta>0$ such that, under the conditions  
\begin{equation}\label{ex(4-2-2-condit-1)}
\quad \left| {x - \overline x } \right| \le \delta, \quad \left\| {h - \overline h } \right\|_\infty   \le \delta, \quad \left( {\overline x ,\overline h } \right) \in G  \times L^\infty  \left( {0,t_1 ;W^{1,\infty } \left( G  \right)} \right)^p, 
\end{equation}
we have
\begin{equation}\label{ex(4-2-2-stima-reg-z-1)}
\left| {z\left( {t,x;h} \right) - z\left( {t,\overline x ;\overline h } \right)} \right| \le 
C_6 \left( {\left| {x - \overline x } \right| + \left\| {h - \overline h } \right\|_\infty  } \right),
\end{equation}
where
\begin{equation}\label{ex(4-2-2-def-C9)}
C_6  =  C\left[ {1 + \left\| {b} \right\|_{L^\infty  (0,t_1;h,\overline h  )}^2  + \left\| {D_{\left( {x,y} \right)} b} \right\|_{L^1 (0,t_1;h, \overline h)}^2  } \right] \times 
\end{equation}
$$
\times \Big[ 1+L^{*2}+\left\| z^* \right\|_{L^\infty(\Gamma _ -  )}^2 
+\left\| {c} \right\|_{L^\infty  (0,t_1;h,\overline h  )}^2+\left\| {a} \right\|_{L^\infty  (0,t_1;h,\overline h  )}^2 +
$$
$$
 {  + \Lambda \left( {D _x h,D _x \overline h } \right)^2
  } {  \big( {\left\| {\nabla _{\left( {x,y} \right)} c} \right\|_{L^1 (0,t_1;h,\overline h )}^2  + \left\| {\nabla _{\left( {x,y} \right)} a} \right\|_{L^1 (0,t_1;h, \overline h)}^2 } \big)} \Big] \times
$$
$$
\times \exp \Big\{ {C \Big[ {\left\| {b} \right\|_{L^\infty  (0,t_1;h,\overline h)} }   + \Lambda \left( {D_x h,D_x \overline h } \right)
  }  
  { {\left\| {D_{\left( {x,y} \right)} b} \right\|_{L^1 (0,t_1;h,\overline h )}  + \left\| {c} \right\|_{L^\infty  (0,t_1;h,\overline h  )} } \Big]} \Big\}
$$
and $C$ is a constant depending on  $B_d$,  $\left\| {D _{\left( {t,x'} \right)} b_d } \right\|_{L_{x_d }^1 (0,1;L_{(t,x')}^\infty  (S_{t_1}' ) )} $.

Hence, $z   \in L_t^\infty  (0,t_1 ;W_{(x,y_{loc} )}^{1,\infty } (G \times \mathbb{R}^p )).$

}

\medskip

{\sc {\bf Proof}}. \

To prove \eqref{ex(4-2-2-stima-reg-z-1)}, it is necessary to obtain preliminary estimates. We first deduce 
\begin{equation}\label{ex(4-2-2-stima-con-i-delta)}
\left| {z\left( {t,x;h} \right) - z\left( {t,\overline x ;\overline h } \right)} \right| \le \tilde \Delta^*  + \tilde \Delta_a + \tilde \Delta_c,
\end{equation}
where
\begin{equation}\label{ex(4-2-2-def-delta-z-*-bis)}
\tilde \Delta^* = \left| {z^* \left( {\tau _ -  \left( {t,x,h} \right),{\rm X}\left( {\tau _ -  \left( {t,x,h} \right);t,x,h} \right)} \right) - z^* \left( {\tau _ -  \left( {t,\overline x ,\overline h } \right),{\rm X}\left( {\tau _ -  \left( {t,\overline x ,\overline h } \right);t,\overline x ,\overline h } \right)} \right)} \right|,  
\end{equation}
\begin{equation}\label{ex(4-2-2-def-delta-a-bis)}
\tilde \Delta_a   = \Big | \int\limits_{\tau _ -  \left( {t,x,h} \right)}^t {  a\left( {s,{\rm X}\left( {s;t,x,h} \right),h\left( {s,{\rm X}\left( {s;t,x,h} \right)} \right)} \right)}ds - 
\end{equation}
$$
+\int\limits_{\tau _ -  \left( {t,\overline x,\overline h} \right)}^t {  a\left( {s,{\rm X}\left( {s;t,\overline x,\overline h} \right),\overline h\left( {s,{\rm X}\left( {s;t,\overline x,\overline h} \right)} \right)} \right)}ds\Big |  
$$
and 
\begin{equation}\label{ex(4-2-2-def-delta-c-bis)}
\tilde \Delta_c   = \Big | \int\limits_{\tau _ -  \left( {t,x,h} \right)}^t {  c\left( {s,{\rm X}\left( {s;t,x,h} \right),h\left( {s,{\rm X}\left( {s;t,x,h} \right)} \right)} \right)z\left( {s,{\rm X}\left( {s;t,x,h} \right);h} \right)}ds - 
\end{equation}
$$
+\int\limits_{\tau _ -  \left( {t,\overline x,\overline h} \right)}^t {  c\left( {s,{\rm X}\left( {s;t,\overline x,\overline h} \right),\overline h\left( {s,{\rm X}\left( {s;t,\overline x,\overline h} \right)} \right)} \right)z\left( {s,{\rm X}\left( {s;t,\overline x,\overline h} \right);\overline h} \right)}ds\Big |.  
$$
Proceeding as for the estimates of $\Delta^*$, $\Delta_a$ and $\Delta_c$  (see Lemma 6.1), we deduce there exist  a positive constant $\tilde C$ indipendent from $t,x,h$ and  $\delta(t,x,h)>0$ such that if \eqref{ex(4-2-2-condit-1)} is verified then we have 
\begin{equation}\label{ex(4-2-2-magg-delta-z-*-bis)}
\tilde \Delta^* \le L^* C_{3} \left( {\left| {x  - \overline x } \right| + \left\| {h - \overline h } \right\|_\infty  } \right),  
\end{equation}
\begin{equation}\label{ex(4-2-2-magg-delta-a-bis)}
\tilde \Delta_a  \le \tilde C \left( {\Lambda \left( {D_x h,D_x h} \right)C_2 \left\| {\nabla _{(x,y)} a} \right\|_{L^1 (0,t_1 ;h, \overline h)}  + C_3\left\| a \right\|_{L^\infty  (0,t_1 ;h,\overline h)} } \right)\left( {\left| {x  - \overline x } \right| + \left\| {h - \overline h } \right\|_\infty  } \right),
\end{equation}
\begin{equation}\label{ex(4-2-2-magg-delta-c-bis)}
\tilde \Delta_c  \le  \tilde C \left[ {\left\| {z^* } \right\|_{L^\infty  \left( {\Gamma _ -  } \right)}  + \left\| a \right\|_{L^1 \left( {0,t_1;h, \overline h} \right)} } \right]\exp  {\left\| c \right\|_{L^1 \left( {0,t_1;h, \overline h} \right)} } \times
\end{equation}
$$
\left( {\Lambda \left( {D_x h,D_x h} \right)C_2 \left\| {\nabla _{(x,y)} c} \right\|_{L^1 (0,t_1 ;h, \overline h)}  + C_3\left\| c \right\|_{L^\infty  (0,t_1 ;h,\overline h)} } \right)\left( {\left| {x  - \overline x } \right| + \left\| {h - \overline h } \right\|_\infty  } \right) +
$$
$$
+\int\limits_{\tau _ -  \left( {t, x , h } \right)\vee \tau _ -  \left( {t,\overline x ,\overline h } \right)}^t {\left\| {c\left( {s, \cdot } \right)} \right\|_{L^\infty  ({G};h,\overline h  )}   } \left| {z\left( {s,{\rm X}\left( {s;t,x,h} \right);h} \right) - z\left( {s,{\rm X}\left( {s;t,\overline x ,\overline h } \right);\overline h } \right)} \right|ds  ,
$$
Then, thanks to \eqref{ex(4-2-2-magg-delta-z-*-bis)}, \eqref{ex(4-2-2-magg-delta-a-bis)} and \eqref{ex(4-2-2-magg-delta-c-bis)} we obtain an estimate for $\left| {z\left( {t,x;h} \right) - z\left( {t,\overline x ;\overline h } \right)} \right|$ (see \eqref{ex(4-2-2-stima-con-i-delta)}); afterwards, applying  Gronwall's lemma  we deduce \eqref{ex(4-2-2-stima-reg-z-1)}. 


\  $ \square$ 
\medskip
\medskip

It is important to observe that, assuming the same hypotheses of Lemma 6.1, we can deduce a stronger regularity result for $z$ respect to the one obtained by Lemma 6.2. Indeed, we have  

{\bf Corollary 6.3.}   \ {\it Let us assume the hyphoteses of Lemma 6.1.  

For every $\left( {t,x,h} \right) \in S_{t_1}  \times L^\infty  \left( {0,t_1 ;W^{1,\infty } \left( G  \right)} \right)^p$ there exists $\delta>0$ such that, under the conditions  
\begin{equation}\label{ex(4-2-2-condit-1)}
\quad \left| {t - \overline t } \right| \le \delta, \quad \left| {x - \overline x } \right| \le \delta, \quad \left\| {h - \overline h } \right\|_\infty   \le \delta, \quad \left( {\overline t,\overline x ,\overline h } \right) \in S_{t_1}  \times L^\infty  \left( {0,t_1 ;W^{1,\infty } \left( G  \right)} \right)^p, 
\end{equation}
we have
\begin{equation}\label{ex(4-2-2-stima-reg-z-2)}
\left| {z\left( {t,x;h} \right) - z\left( {\overline t,\overline x ;\overline h } \right)} \right| \le 
C_{7} \left( {\left| {t - \overline t } \right| +\left| {x - \overline x } \right| + \left\| {h - \overline h } \right\|_\infty  } \right),
\end{equation}
where
\begin{equation}\label{ex(4-2-2-def-c-10)}
C_7  =  C\Big[ 1+L^* \Big]\left[ {1 + \left\| {b} \right\|_{L^\infty  (0,t_1;h,\overline h  )}^4  + \left\| {D_{\left( {x,y} \right)} b} \right\|_{L^1 (0,t_1;h, \overline h)}^4  } \right] \times 
\end{equation}
$$
\times \Big[ 1+L^{*2}+\left\| z^* \right\|_{L^\infty(\Gamma _ -  )}^2 
+\left\| {c} \right\|_{L^\infty  (0,t_1;h,\overline h  )}^2+\left\| {a} \right\|_{L^\infty  (0,t_1;h,\overline h  )}^2 +
$$
$$
 {  + \Lambda \left( {D _x h,D _x \overline h } \right)^2
  } {  \big( {\left\| {\nabla _{\left( {x,y} \right)} c} \right\|_{L^1 (0,t_1;h,\overline h )}^2  + \left\| {\nabla _{\left( {x,y} \right)} a} \right\|_{L^1 (0,t_1;h, \overline h)}^2 } \big)} \Big] \times
$$
$$
\times \exp \Big\{ {C \Big[ {\left\| {b} \right\|_{L^\infty  (0,t_1;h,\overline h)} }   + \Lambda \left( {D_x h,D_x \overline h } \right)
  }  
  { {\left\| {D_{\left( {x,y} \right)} b} \right\|_{L^1 (0,t_1;h,\overline h )}  + \left\| {c} \right\|_{L^\infty  (0,t_1;h,\overline h  )} } \Big]} \Big\}
$$
and $C$ is a constant depending on $t_1$,  $B_d$,  $\left\| {D _{\left( {t,x'} \right)} b_d } \right\|_{L_{x_d }^1 (0,1;L_{(t,x')}^\infty  (S_{t_1}' ) )} $.

Therefore, we deduce that  $z   \in W_{(t,x,y_{loc} )}^{1,\infty } (S_{t_1 }  \times \mathbb{R}^p )
$.
}

\medskip

{\sc {\bf Proof}}. \

To show \eqref{ex(4-2-2-stima-reg-z-2)}, we observe that there exists a positive number $\delta$ such that if $\left| {t - \overline t } \right| \le \delta$, $\left| {x - \overline x } \right| \le \delta$ and $\left\| {h - \overline h } \right\|_\infty   \le \delta 
$, then it follows the estimate 
\begin{equation}\label{ex(4-2-2-new-diff-z)}
\left| {z(t,x;h) - z(\overline t ,\overline x ;\overline h )} \right| \le 
\end{equation}
$$
 \le \left| {z^* (\tau _ -  (t,x,h),{\rm X}(\tau _ -  (t,x,h);t,x,h)) - z^* (\tau _ -  (\overline t ,\overline x ,\overline h ),{\rm X}(\tau _ -  (\overline t ,\overline x ,\overline h );\overline t ,\overline x ,\overline h ))} \right| + 
$$
$$
  \int\limits_{I(t,\overline t ) \cup I(\tau _ -  (t,x,h),\tau _ -  (\overline t ,\overline x ,\overline h ))} {\left[ {\left\| c \right\|_{L^\infty  (0,t_1 ;h,\overline h )} (\left\| {z( \cdot ;h)} \right\|_{L^\infty  (S_{t_1 } )}  \vee \left\| {z( \cdot ;\overline h )} \right\|_{L^\infty  (S_{t_1 } )} ) + \left\| a \right\|_{L^\infty  (0,t_1 ;h,\overline h )} } \right]ds  } 
$$
$$
 + \int\limits_{\tau _ -  (t,x,h) \vee \tau _ -  (\overline t ,\overline x ,\overline h )}^{t \wedge \overline t } {\left[ {\left\| {z( \cdot ;h)} \right\|_{L^\infty  (S_{t_1 } )} \left\| {D_{(x,y)} c(s, \cdot )} \right\|_{L^\infty  (G;h,\overline h )}  + \left\| {D_{(x,y)} a(s, \cdot )} \right\|_{L^\infty  (G;h,\overline h )} } \right]}  \times 
$$
$$
 \times \Big[ {\Lambda (D_x h,D_x \overline h )\left| {{\rm X}(s;t,x,h) - {\rm X}(s;\overline t ,\overline x ,\overline h )} \right| + \left\| {h - \overline h } \right\|_\infty  } \Big]ds + 
$$
$$
 + \int\limits_{\tau _ -  (t,x,h) \vee \tau _ -  (\overline t ,\overline x ,\overline h )}^{t \wedge \overline t } {\left\| {\overline c } \right\|_{L^\infty  (0,t_1 ;\overline h )} \left| {z\left( {s,{\rm X}(s;t,x,h);h} \right) - z\left( {s,{\rm X}(s;\overline t ,\overline x ,\overline h );\overline h } \right)} \right|} ds.
$$
Applying \eqref{ex(3-1-norma-X-second-bis)}, \eqref{ex(3-1-norma-tau-first)} and \eqref{ex(4-2-2-stima-reg-z-1)} to \eqref{ex(4-2-2-new-diff-z)}, we deduce \eqref{ex(4-2-2-stima-reg-z-2)}.

\  $ \square$ 
\medskip
\medskip

{\bf Remark 6.3.0.} The estimate \eqref{ex(4-2-2-stima-reg-z-1)} is very useful to study the IBVP for our quasilinear hyperbolic system (see the proof in Section 7). In particular, the key point is the expression of the constant $C_6$, where the jacobian matrices of $h$ and $\overline h$ (parametric vector functions) are multiplied by norms in $L^1$ respect to time $t_1$; therefore these terms can be chosen small if $t_1$ is assumed sufficiently small.     

\medskip 
We conclude this section with a result of  continuous dependence of the generalized solution on the vector field $b$ and the parameter function $h$.

{\bf Lemma 6.4.}   \ {\it Assume $z^{*(k)} \in Lip^{unif}_{loc} \left( {\Gamma _ -  } \right)
,$ $ c^{(k)}$, $a^{(k)}$ 
$\in L^\infty_t \big( 
{0,t_1;L_{(x,y_{loc} )}^{\infty } \left( {G  \times \mathbb{R}^p} \right)} \big) \cap L^1_t \big(
{0,t_1;W_{(x,y_{loc} )}^{1,\infty } \left( {G  \times \mathbb{R}^p} \right)} \big)$  and let $(t,x,h^{(k)},b^{(k)}) \in S_{t_1} \times  L^\infty  \left( {0,t_1 ;W^{1,\infty } \left( G  \right)} \right)^p
\times L_{x_d }^{-} (S_{t_1 } \times \mathbb{R}^p )
$ with $k=1,2$. If $z^{(k)}$ is the solution of the following integral equation
\begin{equation}\label{ex(4-2-2-int-eq-bar)}
 z^{(k)}(t,x; h^{(k)}) = z^{*(k)} (\tau^{(k)} _ -  (t,x,h^{(k)}),{\rm X}^{(k)}(\tau^{(k)} _ -  ( t,x,h^{(k)});t,x,h^{(k)})) +  
\end{equation}
$$
+\int\limits_{\tau^{(k)} _ -  \left( {t,x,h^{(k)}} \right)}^t {\left[ -{c^{(k)}\left( {s,{\rm X}^{(k)}\left( {s;t,x,h^{(k)}} \right),h^{(k)}\left( {s,{\rm X}^{(k)}\left( {s;t,x,h^{(k)}} \right)} \right)} \right)} \right.}  \times 
$$
$$
 \times z^{(k)}\left( {s,{\rm X}^{(k)}\left( {s;t,x,h^{(k)}} \right); h^{(k)}} \right) + \left. {a^{(k)}\left( {s,{\rm X}^{(k)}\left( {s;t,x,h^{(k)}} \right),h^{(k)}\left( {s,{\rm X}^{(k)}\left( {s;t,x, h^{(k)}} \right)} \right)} \right)} \right]ds,
$$
where ${\rm X}^{(k)}(\cdot;t,x,h^{(k)}) :[ \tau _ - ^{(k)} \left( {t,x,h^{(k)} } \right),\tau _ + ^{(k)} \left( {t,x,h^{(k)} } \right)] \to G$ is the flow associated to the vector field $b^{(k)}$ with $k=1,2$,
then for every $M>0$ there exists $\delta >0$ depending on $h^{\left( 1 \right)},b^{\left( 1 \right)},M$, such that, under the conditions $\left\| {h^{\left( 2 \right)} - h^{\left( 1 \right)} } \right\|_{\infty} < \delta$, ${\left\| {b^{\left( 2 \right)} - b^{\left( 1 \right)} } \right\|_{L^1 (0,t_1;h^{\left( 1 \right)},h^{\left( 2 \right)})} }$
$< \delta$ and $\left\| {b^{(k)} } \right\|_{L^\infty  (0,t_1 ;h^{(k)} )}  \le M
$ with $k=1,2$, we have
\begin{equation}\label{ex(4-4-2-diff-z-z-bar)}
\left\| {z^{(2)}( \cdot ;h^{(2)}) - z^{(1)} ( \cdot ; h^{(1)} )} \right\|_{L^\infty  (S_{t_1} )}  \le 
\end{equation}
$$
\le C_{9} \Big[ {  \left\| {z^{*(2)}  - z^{*(1)}} \right\|_{L^\infty  (\Gamma _ -  )}  + \left\| {b^{(2)} - b^{(1)} } \right\|_{L^1 (0,t_1;h^{(1)}, h^{(2)})} +\left\| {c^{(2)} - c^{(1)} } \right\|_{L^1 (0,t_1;h^{(1)},h^{(2)} )}+  } 
$$
$$
+\left\| {a^{(2)} -  a^{(1)} } \right\|_{L^1 (0,t_1;h^{(1)},h^{(2)} )}+ \left\| {h^{(2)} - h^{(1)} } \right\|_\infty \Big],  
$$
where
\begin{equation}\label{ex(4-2-2-def-c-11)}
C_9  = C \Big( {1 + \mathop  \bigvee \limits_{k = 1,2} \left\| {b^{(k)} } \right\|_{L^\infty  (0,t_1 ;h^{(k)}  )} } \Big) \Big( 1 + \left\| {D_{(x,y)} b^{(1)} } \right\|_{L^1 (0,t_1 ;h^{(1)} ,h^{(2)} )} 
 \Big) \times 
\end{equation}
$$
\times \mathop  \bigvee \limits_{k = 1,2} \Big[ 1+L^{*(k)2}+\left\| {z^{*(k)}} \right\|_{L ^\infty(\Gamma _ -  )}^2 
+\left\| {c^{(k)}} \right\|_{L^\infty  (0,t_1;h^{(k)} )}^2+\left\| {a^{(k)}} \right\|_{L^\infty  (0,t_1;h^{(k)}  )}^2 +
$$
$$
 {  + \Lambda \left( {D _x h^{(1)},D _x h^{(2)} } \right)^2
  } {  \big( {\left\| {\nabla _{\left( {x,y} \right)} c^{(k)}} \right\|_{L^1 (0,t_1;h^{(1)},h^{(2)} )}^2  + \left\| {\nabla _{\left( {x,y} \right)} a^{(k)}} \right\|_{L^1 (0,t_1;h^{(1)}, h^{(2)})}^2 } \big)} \Big] \times
$$
$$
\times \exp \Big[ C \Big( { \mathop  \bigvee \limits_{k = 1,2} \big[\left\| {b^{(k)} } \right\|_{L^\infty  (0,t_1 ;h^{(k)} )} }+ \left\| {c^{(k)} } \right\|_{L^\infty  (0,t_1 ;h^{(k)} )}  \big] + $$
$$
+ \Lambda \left( {D_x h^{(1)} ,D_x h^{(2)} } \right)\left\| {D_{(x,y)} b^{(1)} } \right\|_{L^1 (0,t_1 ;h^{(1)} ,h^{(2)} )} 
   \Big)\Big]
$$   
where the constant $C$  depends only on  $B_d^{\left( k \right)}$, $\left\|{  b_d^{\left( k \right)} } \right\|_{L_{x_d }^1 (0,1;L_{(t,x')}^\infty  (S_{t_1}'  ))} $,  $\left\| {D _{\left( {t,x'} \right)} b_d^{\left( 1 \right)} } \right\|_{L_{x_d }^1 (0,1;L_{(t,x')}^\infty  (S_{t_1}'  ))} $ and $L^{*(k)}$ is the Lipschitz constant for $z^{*(k)}$   $($with $k=1,2)$.

}

{\sc {\bf Proof}}. \

Taking into account Lemma 5.3, the proof is analogous to the one used  in Lemma 6.2.

\  $ \square$ 
\medskip
\medskip
\medskip
\medskip
\section{- The proof of the main theorem.}

Before to prove the main theorem of this paper, we give a preliminary lemma about the semilinear part of our quasilinear hyperbolic system. 

{\bf Lemma 7.1.}   \ {\it Assume \eqref{ex(1-10-f-g)}, \eqref{ex(1-10-condiz-vi)} and \eqref{ex(1-10-condiz-iniz)} for $f_i$, $v_i$ and $y_{0i}$ $(i=1,...,p)$ respectively.     If $C_W$ is a positive constant, then there exists $0<t_2(C_W) \le t_1$ such that the following system of integral equations 
\begin{equation}\label{ex(7-1-IS-eq-y)}
y_i (t,x) = y_{0i} ({\rm X}_{iY} (0;t,x)) + \int\limits_0^t {f_i (s,{\rm X}_{iY} (s;t,x),y(s,{\rm X}_{iY} (s;t,x),\overline w(s,{\rm X}_{iY} (s;t,x)))ds}, 
\end{equation} 
\begin{equation}\label{ex(7-1-IS-eq-X-Y)}
{\rm X}_{iY} (s;t,x) = x - \int\limits_s^t {v_i (r,{\rm X}_{iY} (r;t,x))dr}, \quad s \in [0,t], 
\end{equation} 
where $i=1,..,p$ and $\left\| {\overline w } \right\|_{L^\infty  (0,t_1 ;W^{1,\infty } (G))}  \le C_W$, admits one and only one solution $y \in {L^\infty  (0,t_2 ;W^{1,\infty } (G))^p}$ which satisfies the following estimates
\begin{equation}\label{ex(7-1-norma-y)}
\left\| y \right\|_{L^\infty  (S_{t_2 } )}  \le \left\| {y_0 } \right\|_{L^\infty  (G)}  + 1,\quad \left\| {D_x y} \right\|_{L^\infty  (S_{t_2 } )}  \le 2pd \left\| {D_x y_0 } \right\|_{L^\infty  (G)}  + 1.
\end{equation} 

Furthermore, if  $\overline w^{(k)}$ satisfies
\begin{equation}\label{ex(7-1-norma-w1-w2)}
\left\| {\overline w^{(k)} } \right\|_{L^\infty  (0,t_1 ;W^{1,\infty } (G))}  \le C_W
\end{equation} 
and $y^{(k)} \in {L^\infty  (0,t_2 ;W^{1,\infty } (G))^p}$ is the solution of \eqref{ex(7-1-IS-eq-y)}-\eqref{ex(7-1-IS-eq-X-Y)} with $\overline w =\overline w^{(k)}$ where $k=1,2$, then we have
\begin{equation}\label{ex(7-1-stima-perturb-y)}
\left\| {y^{(2)}  - y^{(1)} } \right\|_{L^\infty  (S_{t_2 } )}  \le \exp (\left\| {D_y f } \right\|_{L^1 (0,t_2 ;(y^{(1)} ,\overline w ^{(1)} ),(y^{(2)} ,\overline w ^{(2)} ))} ) \times
\end{equation}
$$
\times \left\| {D_w f } \right\|_{L^1 (0,t_2 ;(y^{(1)} ,\overline w ^{(1)} ),(y^{(2)} ,\overline w ^{(2)} ))} \left\| {\overline w^{(2)}  - \overline w^{(1)} } \right\|_{L^\infty  (S_{t_2 } )}. 
$$

}

\medskip

{\sc {\bf Proof}}. \

First we consider the following linearized version of \eqref{ex(7-1-IS-eq-y)}  
\begin{equation}\label{ex(7-1-IS-eq-y-lin)}
\tilde y_i (t,x) = y_{0i} ({\rm X}_{iY} (0;t,x)) + \int\limits_0^t {f_i (s,{\rm X}_{iY} (s;t,x),\overline y(s,{\rm X}_{iY} (s;t,x),\overline w(s,{\rm X}_{iY} (s;t,x)))ds}, 
\end{equation} 
where $i=1,...,p,$ $\overline y \in {L^\infty  (0,t_1 ;W^{1,\infty } (G))^p}$ and it satisfies the following estimates
\begin{equation}\label{ex(7-1-norma-overline-y)}
\left\| \overline y \right\|_{L^\infty  (S_{t_1 } )}  \le \left\| {y_0 } \right\|_{L^\infty  (G)}  + 1,\quad \left\| {D_x \overline y} \right\|_{L^\infty  (S_{t_1 } )}  \le 2pd\left\| {D_x y_0 } \right\|_{L^\infty  (G)}  + 1.
\end{equation} 
It is obvious that the problem \eqref{ex(7-1-IS-eq-y-lin)}-\eqref{ex(7-1-IS-eq-X-Y)} admits one and only one solution $\tilde y \in {L^\infty  (0,t_1 ;W^{1,\infty } (G))^p}$ explicitly defined by \eqref{ex(7-1-IS-eq-y-lin)}. Therefore, we can consider the following operator $L_t :L^\infty  (0,t;W^{1,\infty } (G))^p \to L^\infty  (0,t;W^{1,\infty } (G))^p
$ such that $L_t (\overline y) = \tilde y$ with $0 <t \le t_1$. To show that $L_t$ is a contraction for $t$ sufficiently small, we need some estimates about $\tilde y$. More precisely, from \eqref{ex(7-1-IS-eq-y-lin)}-\eqref{ex(7-1-IS-eq-X-Y)}, we obtain 
\begin{equation}\label{ex(7-1-norma-tilde-y)}
\left\| \tilde y \right\|_{L^\infty  (S_{t } )}  \le \left\| {y_0 } \right\|_{L^\infty  (G)}  + \left\| f \right\|_{L^1 (0,t;\overline y, \overline w)}, \quad \left\| {D_x {\rm X}_Y } \right\|_{L^\infty  (S_t )}  \le pd\exp (\left\| {D_x v} \right\|_{L^1 (0,t ;L^\infty  (G))} ),
\end{equation} 
where $0<t \le t_1$ and ${\rm X}_Y = ({\rm X}_{1Y},...{\rm X}_{pY})$. Now, applying the differential operator $D_x$ to \eqref{ex(7-1-IS-eq-y-lin)} and taking into account \eqref{ex(7-1-norma-tilde-y)}, we obtain the following estimates
\begin{equation}\label{ex(7-1-norma-deriv-tilde-y)}
\left\| {D_x \tilde y} \right\|_{L^\infty  (S_t )}  \le  
(pd)\exp (\left\| {D_x v} \right\|_{L^1 (0,t ;L^\infty  (G))} ) 
(\left\| {D_x y_0 } \right\|_{L^\infty  (G)}  +\left\| {D_x f} \right\|_{L^1 (0,t;\overline y,\overline w )} +
\end{equation}
$$
+\left\| {D_y f} \right\|_{L^1 (0,t;\overline y,\overline w )}\left\| {D_x \overline y } \right\|_{L^\infty  (S_{t } )}  + \left\| {D_w f} \right\|_{L^1 (0,t;\overline y,\overline w )} \left\| {D_x \overline w } \right\|_{L^\infty  (S_{t } )} ), \quad 0<t \le t_1.
$$
Hence, from \eqref{ex(7-1-norma-overline-y)}, \eqref{ex(7-1-norma-tilde-y)} and \eqref{ex(7-1-norma-deriv-tilde-y)}, we deduce that there exists $0<t'_1 \le t_1$ such that 
\begin{equation}\label{ex(7-1-norma-tilde-y-bis)}
\left\| \tilde y \right\|_{L^\infty  (S_{t } )}  \le \left\| {y_0 } \right\|_{L^\infty  (G)}  + 1,\quad \left\| {D_x \tilde y} \right\|_{L^\infty  (S_{t } )}  \le 2pd\left\| {D_x y_0 } \right\|_{L^\infty  (G)}  + 1, \quad 0<t\le t'_1;
\end{equation} 
therefore $L_t(C_t)  \subseteq C_t$ if we define $C_t$ as follows
\begin{equation}\label{ex(7-1-def-C-t)}
C_t  = \left\{ {z \in L^\infty  (0,t;W^{1,\infty } (G))^p|} \right.
\end{equation} 
$$
\left. {\left\| z \right\|_{L^\infty  (S_t )}  \le \left\| {y_0 } \right\|_{L^\infty  (G)}  + 1, \quad \left\| {D_x z} \right\|_{L^\infty  (S_t )}  \le 2pd\left\| {D_x y_0 } \right\|_{L^\infty  (G)}  + 1} \right\}, \quad 0<t \le t'_1.
$$
Hence $C_t$ is a closed of $L^\infty  (0,t;W^{1,\infty } (G))^p$. Thus, we have to show that $L_t:C_t \to C_t$ is a contraction; for this purpose, assuming $\overline y= \overline y ^{(k)} \in C_t$ and denoting by $\tilde y ^{(k)} \in C_t$ the solution of the problem \eqref{ex(7-1-IS-eq-y-lin)}-\eqref{ex(7-1-IS-eq-X-Y)} with $k=1,2$, it is not so hard to obtain the following perturbation estimate
\begin{equation}\label{ex(7-1-new-stima-pert)}
\left\| {\tilde y^{(2)}  - \tilde y^{(1)} } \right\|_{L^\infty  (S_t )}  \le \left\| {D_y f } \right\|_{L^1 (0,t;(\overline y ^{(1)} ,\overline w ),(\overline y ^{(2)} ,\overline w ))} \left\| {\overline y ^{(2)}  - \overline y ^{(1)} } \right\|_{L^\infty  (S_t )}. 
\end{equation}
Therefore there exists $0<t_2 \le t'_1$ such that $L_t$ is a contraction; consequently, thanks to Banach-Caccioppoli's fixed point theorem, $L_t$ admits one and only one fixed point denoted by $y$. Of course $y$ is the solution of the problem \eqref{ex(7-1-IS-eq-y)}-\eqref{ex(7-1-IS-eq-X-Y)}. 

To prove \eqref{ex(7-1-stima-perturb-y)}, it is enough to consider the integral expressions for $y^{(1)}$ and  $y^{(2)}$, to make the difference between them and finally to apply  Gronwall's lemma.

\  $ \square$ 
\medskip
\medskip

{\bf Remark 7.1.0.} Of course, in Lemma 7.1, the solution $y$ lies in $W^{1,\infty } (S_{t_2 } )^p$. However, in the proof of the main theorem, first, we prove that the solution $(y,w)$ of our IBVP belongs to ${L^\infty  (0,t^* ;W^{1,\infty } (G))^{p+q}} $; afterwards, it is straightforward to deduce  $(y,w) \in W^{1,\infty } (S_{t^* } )^{p+q}$.

Now we are ready to prove Theorem 3.1.
\medskip
\medskip

{\sc {\bf Proof of the main theorem}}. \

Let us define the following vectors  
\begin{equation}\label{ex(7-10-def-alfa)}
\alpha  = \frac{1}{p}(\left\| {y_0 } \right\|_{L^\infty  (G)}  + 1)\sum\limits_{i = 1}^p {e_i^{(p)} } 
, \quad
\beta  = \frac{1}{q}(\left\| {w^* } \right\|_{L^\infty  (\Gamma _ -  )}  + 1)\sum\limits_{j = 1}^q {e_j^{(q)} } 
\end{equation} 
where $\{ e_1^{(p)}  , ... , e_p^{(p)}\}$, $\{ e_1^{(q)}  , ... , e_q^{(q)}\}$  are the canonical bases of $\mathbb{R}^p$ and $\mathbb{R}^q$ respectively; of course, we have 
\begin{equation}\label{ex(7-10-norma-alfa)}
\left\| \alpha  \right\|_{L^\infty  (S_{t_1 } )}  = \left\| {y_0 } \right\|_{L^\infty  (G)}  + 1,\quad \left\| \beta  \right\|_{L^\infty  (S_{t_1 } )}  = \left\| {w^* } \right\|_{L^\infty  (\Gamma _ -  )}  + 1.
\end{equation}
Moreover, in relation to the vector fields $u_j \in L_{x_d }^ -  ( S_{t_1 })$, we introduce the positive constants $B_{jd}$ as follows
\begin{equation}\label{ex(7-10-def-Bjd)}
u_{jd}(t,x) \le -B_{jd} \quad \mbox{ a.e. } t \in (0,t_1), \quad \forall x \in G.
\end{equation} 
Afterwards, we define the following closed subset of $L^\infty  (0,t;W^{1,\infty } (G))^q$
\begin{equation}\label{ex(9-10-def-C-t)}
K_t  = \left\{ {z \in L^\infty  (0,t;W^{1,\infty } (G))^q|} \right.
\end{equation} 
$$
\left\| z \right\|_{L^\infty  (S_t )}  \le \left\| {w^* } \right\|_{L^\infty  (\Gamma _ -  )}  + 1, \quad \left\| {D_x z} \right\|_{L^\infty  (S_t )}  \le   C\left( {2 + \left\| u \right\|_{L^\infty  (0,t;\alpha )}^2 } \right) \times 
$$
$$
\left. { \times \left( {2 + L^{*2}+\left\| {w^* } \right\|_{L^\infty(\Gamma _ -  )}^2  + \left\| g \right\|_{L^\infty  (0,t;\alpha ,\beta )}^2 } \right)\exp  C\left( {1 + \left\| u \right\|_{L^\infty  (0,t;\alpha )} } \right)
} \right\}, \quad 0<t \le t_1;
$$
where $ C$ is the same constant that appears in the definition of $C_6$ (see \eqref{ex(4-2-2-def-C9)}) but in this case it depends on $\sum\limits_{j = 1}^q {B_{jd} }$,  $\sum\limits_{j = 1}^q \left\| {D _{\left( {t,x'} \right)} u_{jd} } \right\|_{L_{x_d }^1 ( 0,1;L_{( t,x')}^\infty  ( S_{t_1}' ))} $; moreover, $L^*$ is the lipschitz constant for $w^*$.
Thanks to Lemma 7.1, there exists $0<t_2 \le t_1$ such that Cauchy's problem \eqref{ex(7-1-IS-eq-y)}-\eqref{ex(7-1-IS-eq-X-Y)} admits one and only one solution $y \in {L^\infty  (0,t_2 ;W^{1,\infty } (G))^p}$ for every $\overline w \in K_{t_1}$. Now, we consider the operator $H_t:K_t \to L^\infty  (0,t;W^{1,\infty } (G))^q$ such that $H_t(\overline w)= \tilde w $ where $0<t \le t_2$, $\tilde w$ is defined as follows  
\begin{equation}\label{ex(7-10-IS-eq-w)}
\tilde w_j (t,x) = w_j^* (\tau _{j - } (t,x,y),{\rm X}_{jW} (\tau _{j - } (t,x,y);t,x,y)) + 
\end{equation} 
$$
 + \int\limits_{\tau _{j - } (t,x,y)}^t {g_j (s,{\rm X}_{jW} (s;t,x,y),y(s,{\rm X}_{jW} (s;t,x,y)),\overline w(s,{\rm X}_{jW} (s;t,x,y)))ds} ,
$$
\begin{equation}\label{ex(7-10-IS-eq-X-W)}
{\rm X}_{jW} (s;t,x,y) = x - \int\limits_s^t {u_j (r,{\rm X}_{jW} (r;t,x,y),y(r,{\rm X}_{jW} (r;t,x,y)))dr},\quad s \in [\tau _{j - } (t,x,y),t], 
\end{equation} 
$\overline w \in K_t $ and $j=1,...,q $. 

With reference to Lemma 6.2, fixed $\overline w \in K_t$ and $j$, we put
\begin{equation}\label{ex(7-10-def-funzioni)}
z= \tilde w_j, \quad z^*=w^*_j, \quad c=0, \quad h=\overline h =y,\quad a(\cdot,h(\cdot))=g_j (\cdot,y(\cdot),\overline w(\cdot)), \quad b=u_j ;
\end{equation}
therefore the estimates \eqref{ex(4-2-2-stima-z)} and \eqref{ex(4-2-2-stima-reg-z-1)} added from $j=1$ to $j=q$,  give us the following bounds for $\tilde w$
\begin{equation}\label{ex(7-20-stima-tilde-w)}
\left\| {\tilde w} \right\|_{L^\infty  \left( {S_{t} } \right)}  \le \left\| {w^* } \right\|_{L^\infty  \left( {\Gamma _ -  } \right)}  + \left\| g \right\|_{L^1 \left( {0,t;y, \tilde w} \right)} ,
\end{equation}
\begin{equation}\label{ex(4-2-2-stima-deriv-tilde-w)}
\left\| {D_x \tilde w} \right\|_{L^\infty  \left( {S_{t} } \right)} \le   C\left( {1 + \left\| {u} \right\|_{L^\infty  \left( {0,t;y  } \right)}^2  + \left\| {D_{\left( {x,y} \right)} u} \right\|_{L^1 \left( {0,t;y} \right)}^2  } \right) \Big[ 1+L^{*2}+\left\| w^* \right\|_{L^\infty(\Gamma _ -  )}^2 
+ 
\end{equation}
$$
 {+\left\| {g} \right\|_{L^\infty  \left( {0,t;y, \overline w } \right)}^2  + \Lambda \left( {D _x y,D _x y } \right)^2
  } \Big( \left\| {D _{\left( {x,y} \right)} g} \right\|_{L^1 \left( {0,t;y, \overline w} \right)}+ \left\| {D _{w} g} \right\|_{L^1 \left( {0,t;y, \overline w} \right)} \left\| {D_x \overline w} \right\|_{L^\infty  \left( {S_{t} } \right)}\Big)^2 \Big] \times$$
$$
\times \exp \Big[ { C\Big( {\left\| {u} \right\|_{L^\infty  \left( {0,t;y} \right)} }   + \Lambda \left( {D_x y,D_x  y } \right)
  }  
  { {\left\| {D_{\left( {x,y} \right)} u} \right\|_{L^1 \left( {0,t;y} \right)}  } \Big)} \Big].
$$
From these estimates we immediately deduce there exists $0<t'_2\le t_2$ such that if $0 < t \le t'_2$ then we have $H_t (K_t)  \subseteq  K_t$. Hence, to prove the existence and uniqueness of the solution for the system of integral equations \eqref{ex(1-1-IS-eq-y)}-\eqref{ex(1-1-IS-eq-X-W)}, we must prove that $H_t$ is a contraction (for small values of $t$). For this purpose, let $\overline w^{(k)}$ be given in $K_t$ with $0<t \le {t'_2}$ and  $k=1,2$; moreover, we indicate the solution of \eqref{ex(7-1-IS-eq-y)}-\eqref{ex(7-1-IS-eq-X-Y)}, with $\overline w =\overline w^{(k)}$, by $y^{(k)} \in {L^\infty  (0,t ;W^{1,\infty } (G))^p}$. Afterwards, we define $\tilde w^{(k)}=H_t (\overline w^{(k)})$. Now, with reference to Lemma 6.4, we put  
\begin{equation}\label{ex(7-20-def-funzioni)}
z^{(k)} = \tilde w_j^{(k)}, \quad z^{*(k)}=w_j^*, \quad c^{(k)} =0,\quad h^{(k)} =  y^{(k)},  \end{equation}
$$
a^{(k)}(\cdot, h^{(k)}(\cdot)) = g_j(\cdot, y^{(k)}(\cdot),\overline w^{(k)}(\cdot)),\quad b^{(k)} = u_j. 
$$
Furthermore, thanks to the estimate \eqref{ex(7-1-stima-perturb-y)}, there exists $0 < t_3 \le t'_2$ such that
\begin{equation}\label{ex(7-20-diff-y)}
\left\| {y^{(2)}  - y^{(1)} } \right\|_{L^\infty  (S_{t_3 } )}  < \delta 
\end{equation}
where $\delta$ is defined as in Lemma 6.4;
therefore the estimate \eqref{ex(4-4-2-diff-z-z-bar)}  added from $j=1$ to $j=q$,  give us the following perturbation estimate
\begin{equation}\label{ex(7-20-stima-perturb-tildew)}
\left\| {\tilde w^{(2)}  - \tilde w^{(1)} } \right\|_{L^\infty  (S_t )}  \le C\left\{ {\left\| {g ( \cdot ,\overline w ^{(2)} ) - g ( \cdot ,\overline w ^{(1)} )} \right\|_{L^1 (0,t;y^{(1)} ,y^{(2)} )}  + \left\| {y^{(2)}  - y^{(1)} } \right\|_{L^\infty  (S_t )} } \right\}
\end{equation}
where $0<t \le t_3$ and the constant $C$ does not depend from $\overline w^{(k)}$ with $k=1,2$. Now, using \eqref{ex(7-1-stima-perturb-y)} and simple manipulations we deduce
\begin{equation}\label{ex(7-20-stima-perturb-tildew)}
\left\| {\tilde w^{(2)}  - \tilde w^{(1)} } \right\|_{L^\infty  (S_t )}  \le C \left[ {\left\| {D_w g} \right\|_{L^1 (0,t;(y^{(1)} ,\overline w ^{(1)} ),(y^{(2)} ,\overline w ^{(2)} ))}  +  } \right.
\end{equation}
$$
{\left. {+\exp (\left\| {D_y f } \right\|_{L^1 (0,t ;(y^{(1)} ,\overline w ^{(1)} ),(y^{(2)} ,\overline w ^{(2)} ))} )  \left\| {D_w f } \right\|_{L^1 (0,t ;(y^{(1)} ,\overline w ^{(1)} ),(y^{(2)} ,\overline w ^{(2)} ))} 
} \right]} \times
$$
$$
\times \left\| {\overline w^{(2)}  - \overline w^{(1)} } \right\|_{L^\infty  (S_{t } )}, \quad 0<t\le t_3.
$$
Hence there exists $0<t^* \le t_3$ such that $\Lambda_{t^*}$ is a contraction, therefore, thanks to Banach-Caccioppoli's fixed point theorem, this operator admits one and only one fixed point $w$. Therefore we have showed that $(y,w)$, where $y$ is the solution of \eqref{ex(7-1-IS-eq-y)}-\eqref{ex(7-1-IS-eq-X-Y)} with $\overline w=w$, is the generalized solution of  \eqref{ex(1-1-IBVP-eq-y)}-\eqref{ex(1-1-dato-frontiera-w)}. 

Afterwards we study the continuous dependence of generalized solutions on data for the problem \eqref{ex(1-1-IBVP-eq-y)}-\eqref{ex(1-1-dato-frontiera-w)}. For this purpose, we assume \eqref{ex(1-1-IBVP-eq-y)}-\eqref{ex(1-1-dato-frontiera-w)}:
\begin{equation}\label{ex(7-30-f-g)}
f_i ^{(k)}, g_j ^{(k)}  \in 
L^\infty_t (0,t_1 ;L_{(x,y_{loc} ,w_{loc} )}^{\infty } (G  \times \mathbb{R}^p  \times \mathbb{R}^q ))
\cap  L^1_t (0,t_1 ;W_{(x,y_{loc} ,w_{loc} )}^{1,\infty } (G  \times \mathbb{R}^p  \times \mathbb{R}^q )),  
\end{equation} 
\begin{equation}\label{ex(7-30-condiz-ui-vi)}
v_i ^{(k)} \in L^1 (0,t_1 ;W^{1,\infty } (G))^d,\quad u_j ^{(k)}  \in L_{x_d }^ -  \left( {S_{t_1 }  \times \mathbb{R}^p }  \right),
\end{equation}
\begin{equation}\label{ex(7-30-condiz-vi)}
v_i ^{(k)} (t,x',0) \cdot e_d = v_i ^{(k)} (t,x',1) \cdot e_d=0 \quad \forall x' \in \mathbb{R}^{d-1},  \mbox{ a.e. } t \in (0,t_1), 
\end{equation}
\begin{equation}\label{ex(7-30-condiz-iniz)}
y_{0i}^{(k)}  \in W^{1,\infty } (G),\quad w_j^{*(k)}  \in Lip^{unif}_{loc}(\Gamma _ -  ), \quad i=1,...,p, \quad j=1,...,q, \quad k=1,2.
\end{equation}
With reference to \eqref{ex(1-1-IS-eq-y)}-\eqref{ex(1-1-IS-eq-X-W)}, we assume $v_i=v_i ^{(k)}$,                $f_i=f_i^{(k)}$, $y_{0i}=y_{0i}^{(k)}$, $u_{j}=u_{j}^{(k)}$, $g_{j}=g_{j}^{(k)}$,  
 $w_{j}^*=w_{j}^{*(k)}$ and we denote by $(y_{i}^{(k)},w_{j}^{(k))})$ the solution of \eqref{ex(1-1-IS-eq-y)}-\eqref{ex(1-1-IS-eq-X-W)} in $(0,\tilde t)$ with $i=1,...,p$, $j=1,...,q$ and $k=1,2$. Subtracting  \eqref{ex(1-1-IS-eq-y)}, \eqref{ex(1-1-IS-eq-X-Y)} with $k=2$ from the same equations with $k=1$, we obtain, after some calculations, the following estimate
\begin{equation}\label{ex(7-30-perturbaz-y)}
\left\| {y^{(2)} (t, \cdot ) - y^{(1)} (t, \cdot )} \right\|_{L^\infty  (G)}  \le  \left\| {y_0^{(2)}  - y_0^{(1)} } \right\|_{L^\infty  (G)}  + \left\| {f^{(2)}  - f^{(1)} } \right\|_{L^1 (0,t;(y^{(1)} ,w^{(1)} ),(y^{(2)} ,w^{(2)} ))}  + 
\end{equation} 
$$
 +\big (\left\| {D_x y_0^{(1)} } \right\|_{L^\infty  (G)}  +\left\| {D_x f^{(1)} } \right\|_{L^1 (0,t;(y^{(1)} ,w^{(1)} ),(y^{(2)} ,w^{(2)} ))}  + \left\| {D_y f^{(1)} } \right\|_{L^1 (0,t;(y^{(1)} ,w^{(1)} ),(y^{(2)} ,w^{(2)} ))} \times
$$
$$
 \times \left\| {D_x y^{(1)} } \right\|_{L^\infty  (S_t )}  +   \left\| {D_w f^{(1)} } \right\|_{L^1 (0,t;(y^{(1)} ,w^{(1)} ),(y^{(2)} ,w^{(2)} ))} \left\| {D_x w^{(1)} } \right\|_{L^\infty  (S_t )} \big)\times 
$$
$$
\times \left\| {v^{(2)}  - v^{(1)} } \right\|_{L^1 (0,t;L^\infty  (G))}\exp \left\| {D_x v^{(1)} } \right\|_{L^1 (0,t;L^\infty  (G))} 
  + 
$$
$$
+ \int\limits_0^t {\left\| {D_y f^{(1)} (s)} \right\|_{L^1 (0,t;(y^{(1)} ,w^{(1)} ),(y^{(2)} ,w^{(2)} ))} } \left\| {y^{(2)} (s, \cdot ) - y^{(1)} (s, \cdot )} \right\|_{L^\infty  (G)} ds + 
$$
$$
 + \int\limits_0^t {\left\| {D_w f^{(1)} (s)} \right\|_{L^1 (0,t;(y^{(1)} ,w^{(1)} ),(y^{(2)} ,w^{(2)} ))} } \left\| {w^{(2)} (s, \cdot ) - w^{(1)} (s, \cdot )} \right\|_{L^\infty  (G)} ds .
$$
Now, with reference to Lemma 6.4, we assume
\begin{equation}\label{ex(7-40-funzioni-k)}
z^{(k)}  = w_j^{(k)} ,\quad z^{*(k)}  = w_j^{*(k)} ,\quad c^{(k)}  = 0,\quad h^{(k)}  = y^{(k)} ,
\end{equation}
$$
a^{(k)}(\cdot, h^{(k)}(\cdot)) = g_j^{(k)}(\cdot, y^{(k)}(\cdot), w^{(k)}(\cdot)),\quad b^{(k)} = u_j^{(k)}; 
$$
therefore, applying \eqref{ex(4-4-2-diff-z-z-bar)} and adding it from $j=1$ to $q$, we obtain 
\begin{equation}\label{ex(7-30-perturbaz-w)}
\left\| {w^{(2)} (t, \cdot ) - w^{(1)} (t, \cdot )} \right\|_{L^\infty  (G)}  \le C \Big(\left\| {w^{*(2)}  - w^{*(1)} } \right\|_{L^\infty  (\Gamma _ -  )}  + \left\| {u^{(2)}  - u^{(1)} } \right\|_{L^1 (0,t;y^{(1)} ,y^{(2)} )}  + 
\end{equation} 
$$
 + \left\| {g^{(2)}  - g^{(1)} } \right\|_{L^1 (0,t;(y^{(1)} ,w^{(1)} ),(y^{(2)} ,w^{(2)} ))}  + \left\| {y^{(2)}  - y^{(1)} } \right\|_{L^\infty  (S_t )} \Big).
$$
Hence, combining \eqref{ex(7-30-perturbaz-y)} with \eqref{ex(7-30-perturbaz-w)}, we deduce the result of continuous dependence.

Finally, remembering Corollary 6.3 and Remark 7.1.0, we deduce that $(y,w) \in W^{1,\infty } (S_{t^* } )^{p+q}$.

\  $ \square$ 
\medskip 
\medskip
\medskip
\medskip
\section{- The well-posedness of an IBVP for the hyperbolic part of an atmospheric model.}

In the paper \cite{[SF]}  we have introduced and studied a model of the phase transitions for $H_2O$ in the three states in the atmosphere. Let us say something about this model. We can consider a spatial domain $\Omega$ in the atmosphere. Inside $\Omega$ there are many particles interacting among them characterized by densities and velocities. More exactly, we consider the density $\rho$ of dry air, the density $\pi$ of water vapour, the density $\sigma$ of $H_2O$ in liquid state, the density $\nu$ of $H_2O$ in solid state, the speed $v$ and the temperature $T$ of the atmospheric gas, the speed $u$ of water droplets and the speed $w$ of ice crystalls. In the model proposed in \cite{[SF]}, these unknown quantities are linked among them by six partial differential equations, of which four of them are hyperbolic type whereas the others are  parabolic type, whereas the speed $u$ of droplets and the speed $w$ of ice crystals are given by simple formulas involving $v$.

In this section, using the main theorem of this paper, we study IBVP for the hyperbolic part of the model in \cite{[SF]} on the strip $ \Omega_M =\left\{ (m,x_1,x_2,x_3) |\quad m,x_1,x_2 \in \mathbb{R}, \quad 0 < x_3 < 1  \right\}$ with given velocities $v,$ $u,$ $w,$ and temperature $T$. Moreover, we assume $v$ to be tangent to the planes $x_3=0$ and $x_3=1$, whereas $u$ and $w$ have negative vertical components; therefore, rain and ice fall from the strip. Hence  we consider the following system of partial differential equations in the unknown functions $\rho$, $\pi$, $\sigma$, $\nu$

\begin{equation}\label{ex(2-1-equaz-in-rho)}  
\partial_t \rho \left( {t, x} \right) + v\left( {t, x} \right) \cdot \nabla _{ x} \rho \left( {t, x} \right) = R^* \left( \rho  \right)\left( {t, x} \right),
\end{equation}
\begin{equation}\label{ex(2-1-equaz-in-pi)}
\partial_t \pi \left( {t, x} \right) + v\left( {t, x} \right) \cdot \nabla _{ x} \pi \left( {t, x} \right) = P^* \left( {\pi ,\sigma ,\nu} \right)\left( {t, x} \right),
\end{equation} 
\begin{equation}\label{ex(2-1-equaz-in-sigma)}
\partial_t \sigma\left( {t,m,x} \right) + ( s_l \left( {m } \right)\left[ {\left( {\pi  - \pi _l \left( T \right)} \right)\left( {t, x} \right)} \right],{u \left( {t,m,x} \right)} ) \cdot \nabla _{(m,x)} \sigma \left( {t,m,x} \right) =   
\end{equation}
$$
=S^* 
 \left( {\pi ,\sigma  ,\nu } \right)\left( {t,m,x} \right),
$$
\begin{equation}\label{ex(2-1-equaz-in-nu)}
\partial_t \nu\left( {t,m,x} \right) + ( s_s \left( {m } \right)\left[ {\left( {\pi  - \pi _s \left( T \right)} \right)\left( {t, x} \right)} \right],{w \left( {t,m,x} \right)} ) \cdot \nabla _{(m,x)} \nu \left( {t,m,x} \right) =   
\end{equation}
$$
=N^* 
 \left( {\pi ,\sigma  ,\nu } \right)\left( {t,m,x} \right),
$$
where $\left( {t,m,x} \right) \in S_{t_1 }  =\left( {0,t_1 } \right) \times \mathbb{R} \times \mathbb{R}^2 \times (0,1) = \left( {0,t_1 } \right) \times \Omega _M=\left( {0,t_1 } \right) \times \mathbb{R} \times \Omega=S'_{t_1 } \times (0,1)  $ and $Q_{t_1 }=\left( {0,t_1 } \right) \times \Omega$; moreover we assume the following conditions 
\begin{equation}\label{ex(2-1-condiz-iniz-rho)}
\rho \left( {0, x} \right) = \rho _0 \left( { x} \right)\quad  x \in  \Omega, 
\end{equation}
\begin{equation}\label{ex(2-1-condiz-iniz-pi)}
\pi \left( {0, x} \right) = \pi _0 \left( { x} \right)\quad  x \in  \Omega, 
\end{equation}
\begin{equation}\label{ex(2-1-condiz-sigma-j)}
\sigma  \left( {t,m,x} \right) = \sigma^* \left( {t,m,x} \right) \quad \left( {t,m,x} \right) \in \Gamma _ -  
,  
\end{equation}
\begin{equation}\label{ex(2-1-condiz-nu-j)}
\nu\left( {t,m,x} \right) = \nu^* \left( {t,m,x} \right) \quad \left( {t,m,x} \right) \in \Gamma _ - 
,  
\end{equation}
where $\Gamma_-$ is the surface carrying data defined as $\left( {\left\{ 0 \right\} \times \Omega_M } \right)$ $ \cup \left( {\left[ {0,t_1} \right] \times \mathbb{R}^3  \times \left\{ 1 \right\}} \right)
$ and $\rho_0$, $\pi_0$, $\sigma^*$, $\nu^*$ are  given functions.

Let us recall the quantities which figure in IBVP \eqref{ex(2-1-equaz-in-rho)}-\eqref{ex(2-1-condiz-nu-j)}. The variable $t$ is the time,  $m$ is the mass of a water droplet or an ice crystal, whereas $x$ is the position vector; in particular, $x_3$ is the height of a generic point identified by $ x$. Moreover, knowing that there are no droplets or ice crystals with smaller mass than $m_a>0$ (see \cite{[SF]}), it will not be restrictive to extend the functions in the partial differential equations of our IBVP to zero for $m \le 0$; this makes it possible to assume $m \in \mathbb{R}$. 

The quantities $\rho_0$, $\pi_0$ are initial densities (for $t=0$), whereas $\sigma^* $, $\nu^*$ are the prescribed densities on $\Gamma_-$. 

Now, we define the functions appearing in the second members of \eqref{ex(2-1-equaz-in-rho)}-\eqref{ex(2-1-equaz-in-nu)}         
\begin{equation}\label{ex(9-1-def-R-*)}
R^* \left( \rho  \right)(t,x) =  -\big[ \left( {\nabla _{ x}  \cdot v} \right)\rho\big](t,x) ;
\end{equation}
\begin{equation}\label{ex(9-1-def-P-*)}
P^* \left( {\pi ,\sigma ,\nu } \right)(t,x) =  - \big[ \left( {\nabla _{ x}  \cdot v} \right)\pi  + P\left( {\pi ,\sigma ,\nu } \right) \big](t,x),
\end{equation}
where $P\left( {\pi ,\sigma ,\nu } \right)$ represents the total amount of water vapour that is transformed into liquid or solid state and we assume that:
\begin{equation}\label{ex(9-1-def-P)}
P\left( {\pi ,\sigma ,\nu } \right)\left( {t,x} \right)=- \big[{\left( {\pi  - \pi _l \left( T \right)} \right)F_l \left( {\sigma } \right)}- {\left( {\pi  - \pi _s \left( T \right)} \right)F_s \left( {\nu } \right)}\big](t,x),
\end{equation}
\begin{equation}\label{ex(9-1-def-F-l)}
F_l \left( {\sigma } \right)\left( {t, x} \right) = \int\limits_0^\infty  {\overline s _l \left( m \right)\sigma \left( {t,m, x} \right)} dm,
\end{equation}
\begin{equation}\label{ex(9-1-def-F-s)}
F_s \left( {\nu } \right)\left( {t, x} \right) = \int\limits_0^\infty  {\overline s _s \left( m \right)\nu \left( {t,m, x} \right)} dm;
\end{equation}
\begin{equation}\label{ex(9-1-def-S-*)}
S^*(\pi,\sigma,\nu)(t,m,x) =  -  \big[\left( {\nabla _{ (x,m)}  \cdot u} \right)\sigma\big](t,m,x)+
\end{equation}
$$   +[S_g(\pi,\sigma) +
S_s(\sigma,\nu) + S_a(\pi,\sigma,\nu) + S_q(\sigma,\nu)](t,m,x),
$$
where
\begin{equation}\label{ex(9-1-def-Sq)}
S_g(\pi,\sigma)(t,m,x) = {\overline s _l(m)\,[\pi(t,x)-\pi_l(T(t,x))]\, \sigma(t,m,x)}
\end{equation}
is the amount of $H_2O$ converted from gas to liquid  that condenses on droplets  with mass $m$,
\begin{equation}\label{ex(9-1-def-S-s)}
S_s(\sigma,\nu)(t,m,x) = [- K_{ls}(m,T)\,\sigma(m) + K_{sl}(m,T)\,\nu(m)](t,x)
\end{equation}
is the amount of droplets with mass $m$ that appears or disappears due to the solidification or fusion,
\begin{equation}\label{ex(9-1-def-S-a)}
S_a(\pi,\sigma,\nu)(t,m,x) =
\end{equation} 
$$= g_a (m) \big[ N^*(t,x)  - 
\int\limits_0^{ + \infty } {n_l (m)\sigma (t,m,x)dm - } \int\limits_0^{ + \infty } {n_s (m)\nu (t,m,x)dm}]^+ 
[\pi(t,x) - \pi_l(T(t,x))]^+  - $$
$$- g_l(m) [ {\pi(t,x)  - \pi_l(T(t,x))} ]^-  \sigma(t,m,x) $$
is the total amount of droplets with mass $m$ which arises or evaporates on the aerosol particles,
\begin{equation}\label{ex(9-1-def-S-q)}
S_q(\sigma,\nu)(t,m,x)=\big[ Q_l(\sigma,\sigma) + J_l(\sigma)\sigma + J_{ls}(\nu)\sigma \big](t,m,x)
\end{equation}
is relative to interactions among droplets and ice crystals and we give the following additional definitions
\begin{equation}\label{ex(9-1-def-J-l)}
J_l (\sigma)(t,m,x) = - m {\int_0^\infty \beta_l (m,m') \sigma (t,m',x)\, dm'},
\end{equation}
\begin{equation}\label{ex(9-1-def-Q-l)}
Q_l ( \sigma, \sigma ')(t,m,x ) \!=\! 
{\frac{m}{2}\! \int_0^m \!\!\beta_l (m'\!,m-m')\sigma (t,m',x) \sigma '(t,m-m',x)dm'\!},
\end{equation}
\begin{equation}\label{ex(9-1-def-J-ls)}
J_{ls} ( \omega)(t,m,x ) =  - m { \int_0^{\infty } Z_{ls}(m'\!,m) \,\omega(t,m',x)\,dm'} , \quad \omega = \nu, \sigma,
\end{equation}
where $\beta_l(m,m')$ is related to the probability 
that droplets of mass $m$ and $m'$ collide and merge,
whereas $Z_{ls}(m,m')$ regards the probability that a droplet of mass $m'$ joins an ice particle of mass $m$ (with instantaneous phase transition from the liquid to the solid state);
\begin{equation}\label{ex(9-1-def-N-*)}
N^*(\pi,\sigma,\nu)(t,m,x) =  -  \big[\left( {\nabla _{ (x,m)}  \cdot w} \right)\nu\big](t,m,x)+
\end{equation}
$$   +[N_g(\pi,\nu) +N_a(\pi,\nu)+
N_s(\sigma,\nu) + N_q(\sigma,\nu)](t,m,x),
$$
where
\begin{equation}\label{ex(9-1-def-N-q)}
N_g(\pi,\nu)(t,m,x) = {\overline s _s(m)\,[\pi(t,x)-\pi_s(T(t,x))]\, \nu(t,m,x)},
\end{equation}
\begin{equation}\label{ex(9-1-def-N-s)}
N_s(\sigma,\nu)(t,m,x) = [ K_{ls}(m,T)\,\sigma(m) - K_{sl}(m,T)\,\nu(m)](t,x),
\end{equation}
\begin{equation}\label{ex(9-1-def-N-a)}
N_a(\pi,\nu)(t,m,x) =  - g_s (m) \big[\pi (t,x) - \pi _s (T(t,x)) \big]^- \nu (t,m,x),
\end{equation}
\begin{equation}\label{ex(9-1-def-N-q)}
N_q(\sigma,\nu)(t,m,x)=\big[ Q_s(\nu,\nu) + J_s(\nu)\nu +Q_{ls}(\nu,\sigma) +J_{ls}(\sigma)\nu \big](t,m,x),
\end{equation}
\begin{equation}\label{ex(9-1-def-J-s)}
J_s (\nu)(t,m,x) = - m {\int_0^\infty \beta_s (m,m') \nu (t,m',x)\, dm'},
\end{equation}
\begin{equation}\label{ex(9-1-def-Q-s)}
Q_s ( \nu, \nu ')(t,m,x ) \!=\! 
{\frac{m}{2}\! \int_0^m \!\!\beta_l (m'\!,m-m')\nu (t,m',x) \nu '(t,m-m',x)dm'\!},
\end{equation}
\begin{equation}\label{ex(9-1-def-Q-ls)}
Q_{ls} ( \nu, \sigma)(t,m,x ) \!=\! 
{\frac{m}{2}\! \int_0^m \!\!Z_{ls} (m'\!,m-m')\sigma (t,m',x) \nu(t,m-m',x)dm'\!}.
\end{equation}
These last quantities are defined in a similar way as already we have seen above about droplets. Moreover, it is possible to find the definitions about the physical quantities such as $s_j$, $\overline s_j$, $\pi_j$, $N^*$, $n_j$, $g_a$, $g_j$, $K_{ls}$ and $K_{sl}$ ($j=l,s$) in \cite{[SF]}, \cite{[AS]}, etc. 

To the  initial and boundary value problem \eqref{ex(2-1-equaz-in-rho)} - \eqref{ex(2-1-condiz-nu-j)} we can associate, by the method of characteristics, the following system of integral equations   
\begin{equation}\label{ex(3-2-def-eq-intg-in-rho)}
\rho \left( {t, x} \right) = \rho_0 ( {\rm X}_\Pi  \left( {0;t, x} \right)) + \int\limits_{0}^t {R^* \left( \rho  \right)} \left( {r,{\rm X}_\Pi  \left( {r;t, x} \right)} \right)dr,
\end{equation}
\begin{equation}\label{ex(3-2-def-eq-intg-in-pi)}
\pi \left( {t, x} \right) =  \pi_0 ( {\rm X}_\Pi  \left( {0;t, x} \right))  + \int\limits_{0}^t {P^* \left( {\pi ,\sigma ,\nu } \right)} \left( {r,{\rm X}_\Pi  \left( {r;t, x} \right)} \right)dr,
\end{equation}
\begin{equation}\label{ex(3-2-def-eq-intg-in-sigma)}
\sigma \left( {t,m,x} \right) =  \sigma^* \left( {\tau _{\Sigma  - } \left( {t,m,x;\pi } \right),{\rm X}_{\Sigma  } \left( {\tau _{\Sigma   - } \left( {t,m,x;\pi } \right);t,m,x;\pi } \right)} \right) + 
\end{equation}
$$
 + \int\limits_{\tau _{\Sigma   - } \left( {t,m,x;\pi } \right)}^t {S^* \left( {\pi ,\sigma ,\nu } \right)} \left( {r,{\rm X}_{\Sigma  } \left( {r;t,m,x;\pi } \right)} \right)dr,
$$
\begin{equation}\label{ex(3-2-def-eq-intg-in-nu)}
\nu \left( {t,m,x} \right) =  \nu^* \left( {\tau _{N  - } \left( {t,m,x;\pi } \right),{\rm X}_{N  } \left( {\tau _{N   - } \left( {t,m,x;\pi } \right);t,m,x;\pi } \right)} \right) + 
\end{equation}
$$
 + \int\limits_{\tau _{N   - } \left( {t,m,x;\pi } \right)}^t {N^* \left( {\pi ,\sigma ,\nu } \right)} \left( {r,{\rm X}_{N  } \left( {r;t,m,x;\pi } \right)} \right)dr,
$$
where $\left( {t,m,x} \right) \in S_{t_1 }$, the fluxes ${\rm X}_\Pi$, ${\rm X}_{\Sigma  }$, ${\rm X}_{N  }$  are defined as follows
\begin{equation}\label{ex(3-2-def-eq-intg-in-X-pi)}
{\rm X}_\Pi  \left( {r;t, x} \right) =  x + \int\limits_t^r v \left( {q,{\rm X}_\Pi  \left( {q;t, x} \right)} \right)dq,\quad r \in \left[ {0,t_1} \right],
\end{equation}
\begin{equation}\label{ex(3-2-def-eq-intg-in-X-sigma)}
{\rm X}_{\Sigma  } \left( {r;t,m,x;\pi } \right) = (m,x) + \int\limits_t^r {\left( {s_l \left( {\pi  - \pi _l \left( T \right)} \right),u } \right)} \left( {q,{\rm X}_{\Sigma  } \left( {q;t,m,x; \pi} \right)} \right)dq,
\end{equation}
$$
r \in \left[ {\tau _{\Sigma   - } \left( {t,m,x;\pi } \right),t} \right], 
$$
\begin{equation}\label{ex(3-2-def-eq-intg-in-X-nu)}
{\rm X}_{N  } \left( {r;t,m,x;\pi } \right) = (m,x) + \int\limits_t^r {\left( {s_s \left( {\pi  - \pi _s \left( T \right)} \right),w } \right)} \left( {q,{\rm X}_{N  } \left( {q;t,m,x; \pi} \right)} \right)dq,
\end{equation}
$$
r \in \left[ {\tau _{N  - } \left( {t,m,x;\pi } \right),t} \right], 
$$
and $\tau _{\Sigma   - } \left( {t,m,x;\pi } \right)$, $\tau _{N   - } \left( {t,m,x;\pi } \right)$ are the minimal time of existence for the solutions of \eqref{ex(3-2-def-eq-intg-in-X-sigma)}, \eqref{ex(3-2-def-eq-intg-in-X-nu)} respectively.

Moreover any solution of the integral equations system \eqref{ex(3-2-def-eq-intg-in-rho)}-\eqref{ex(3-2-def-eq-intg-in-X-nu)} will be considered as a generalized solution for the IBVP \eqref{ex(2-1-equaz-in-rho)}-\eqref{ex(2-1-condiz-nu-j)}. 

We make now assumptions on the functions that appear in \eqref{ex(2-1-equaz-in-rho)}-\eqref{ex(2-1-equaz-in-nu)}. More exactly, we assume the following conditions on $T$, $v$, $u$ and $w$:
\begin{equation}\label{ex(3-2-condiz-temperatura)}
T \in 
L^\infty \left( {0,t_1;L^{\infty } \left( \Omega  \right)} \right)
\cap L^1 \left( {0,t_1;W^{1,\infty } \left( \Omega  \right)} \right), 
\end{equation}
\begin{equation}\label{ex(3-2-ipotes-v)}
v \in L^\infty \left( {0,t_1;L^{\infty } \left( \Omega  \right)} \right)^3 \cap L^1 \left( {0,t_1;W^{1,\infty } \left( \Omega  \right)} \right)^3, 
\end{equation}
\begin{equation}\label{ex(3-2-def-condiz-v)}
v\left( {t,x_1,x_2,0} \right) \cdot e_3  = v\left( {t,x_1,x_2,1} \right) \cdot e_3  = 0,\quad \forall (x_1,x_2) \in \mathbb{R}^2 \quad \mbox{ for almost every } t \in(0,t_1) 
\end{equation}
($\left\{ {e_k |k = 1,2,3} \right\}$ is the canonical frame on $\mathbb{R}^3$),
\begin{equation}\label{ex(3-2-def-condiz-u-j)}
u,w  \in 
L^\infty \left( {0,t_1;L^{\infty } \left( \Omega_M  \right)} \right)^3 \cap L^1 \left( {0,t_1;W^{1,\infty } \left( \Omega_M  \right)} \right)^3, 
\end{equation}
\begin{equation}\label{ex(3-2-def-condiz-u-j-a-bis)}
u_3,w_3  \in  L^1_{x_3} ( 0,1;W^{1,\infty }_{(t,m,x_1,x_2)} ( S'_{t_1}  )), 
\end{equation}
\begin{equation}\label{ex(3-2-u3-v3)}
u_3(t,m,x),w_3(t,m,x) \le -B \quad \forall (m,x) \in \Omega_M \mbox{ a.e. } t \in (0,t_1) ,
\end{equation}
where $B$ is a positive constant $($we say that $u,w \in \tilde L_{x_3 }^ -  \left( {S_{t_1 } } \right)$ for reasons of convenience$)$,
\begin{equation}\label{3-2-def-nabla}
\nabla_{x} \cdot v \in 
L^\infty \left( {0,t_1;L^{\infty } \left( \Omega  \right)} \right)
\cap L^1 \left( {0,t_1;W^{1,\infty } \left( \Omega  \right)} \right),
\end{equation}
\begin{equation}\label{3-2-def-nabla-bis}
\nabla_{(m,x)}  \cdot u,\nabla_{(m,x)}  \cdot w \in 
L^\infty \left( {0,t_1;L^{\infty } \left( \Omega_M  \right)} \right)
\cap  L^1 \left( {0,t_1;W^{1,\infty } \left( \Omega_M  \right)} \right).
\end{equation}

As for initial and boundary conditions, we take
\begin{equation}\label{ex(3-2-def-condiz-iniz-rho-pi)}
\rho _0 ,\pi _0  \in W^{1,\infty } \left( \Omega  \right), \quad \rho _0 ,\pi _0 \geq 0, 
\end{equation}
\begin{equation}\label{ex(3-2-def-condiz-iniz-sigma-*)}
 \sigma^* ,\nu^* \in Lip^{unif}_{loc} \left( {\Gamma _ -  } \right), \quad  \sigma^* , \nu^* \geq 0,
, \quad \sigma^* \left( { \cdot ,m , \cdot } \right) =\nu^* \left( { \cdot ,m , \cdot } \right) =0 \quad \forall m  \notin \left[ {m_a ,M^* } \right],
\end{equation}
where $0<m_a<M^*$.
Moreover, according to the physical model introduced in \cite{[SF]}, we also assume 
\begin{equation}\label{ex(3-3-mt-funz-aus)}
s_j , \overline s_j, g_{a}, g_j ,\pi _j, K_{ji}  \in W^{1,\infty } \left( {\mathbb{R},\mathbb{R}_+ } \right), \quad  \beta _j ,Z_{ls}  \in W^{1,\infty } \left( {\mathbb{R}^2,\mathbb{R}_+   } \right),  
\end{equation}
$$
 n _j  \in L^1 \left( {\mathbb{R},\mathbb{R}_+ } \right),\quad   N^*  \in W^{1,\infty } \left( {(0,t_1) \times \Omega,\mathbb{R}_+  } \right),
$$
$$
supp(g_{0l}),supp(g_j)  \subseteq \left[ {m_a ,M_a } \right], \quad  
 supp(s_j), supp(\overline s_j)  \subseteq \left[ {m_a ,M^* } \right],
$$  
$$ 
\beta _j \left( {m',m''} \right) = Z_{ij} \left( {m',m''} \right) = 0, \quad m' + m'' \ge M^*,
$$
with $\left( {j,i} \right) \in \left\{ {\left( {l,s} \right),\left( {s,l} \right)} \right\}
$  and $0<m_a<M_a<M^*$.



We can now present a theorem of well-posedness about IBVP \eqref{ex(2-1-equaz-in-rho)}-\eqref{ex(2-1-condiz-nu-j)} 
\medskip
\medskip

{\bf Theorem 8.1.}  \ {\it Assume that the hypotheses  \eqref{ex(3-2-condiz-temperatura)}-\eqref{ex(3-3-mt-funz-aus)}   are verified. Then there exists $0 < t^*  \le t_1$ such that  
\eqref{ex(3-2-def-eq-intg-in-rho)}-\eqref{ex(3-2-def-eq-intg-in-nu)} admits one and only one solution $(\rho, \pi, \sigma, \nu) \in $  $W^{1,\infty } \left( {Q_{t^* } } \right)^2  \times W^{1,\infty } \left( {S_{t^* } } \right)^2$; 
 furthermore, $\rho,\pi,\sigma,\nu$ take non-negative values and satisfy the relation 
\begin{equation}\label{ex(3-mt-supp-sigma-nu)}
\sigma \left( { \cdot ,m, \cdot} \right) = \nu \left( { \cdot ,m, \cdot} \right) = 0 \quad \forall m \notin \left[ {m_a ,M^* } \right].
\end{equation} 
The function $(\rho,\pi,\sigma,\nu)$ is also said to be the generalized solution for the equations
\eqref{ex(2-1-equaz-in-rho)}-\eqref{ex(2-1-equaz-in-nu)}  with the initial conditions \eqref{ex(2-1-condiz-iniz-rho)}-\eqref{ex(2-1-condiz-iniz-pi)} and the boundary conditions \eqref{ex(2-1-condiz-sigma-j)}-\eqref{ex(2-1-condiz-nu-j)}.

Moreover, for any sufficiently small t, the mapping $\big( \rho_0 ,\pi_0 ,\sigma^* ,\nu^*,T,v,u,w, \nabla_x \cdot v,$ $ \nabla_{(x,m) } \cdot u, \nabla_{(x,m) } \cdot w \big) \in
$ $W^{1,\infty } \left( {\Omega } \right)^2  \times Lip^{unif}_{loc} \left( {\Gamma_- } \right)^2  \times 
\big[L^\infty \left( {0,t;L^{\infty } \left( \Omega  \right)} \right)
\cap L^1 \left( {0,t;W^{1,\infty } \left( \Omega  \right)} \right) \big]^4
\times \tilde L_{x_3 }^ -  \left( {S_{t } } \right)^2 \times  
\big[L^\infty \left( {0,t;L^{\infty } \left( \Omega  \right)} \right)
\cap L^1 \left( {0,t;W^{1,\infty } \left( \Omega  \right)} \right) \big]
\times 
\big[L^\infty \left( {0,t;L^{\infty } \left( \Omega_M  \right)} \right)
\cap L^1 \left( {0,t;W^{1,\infty } \left( \Omega_M  \right)} \right) \big]^2$ $ \to \left( {\rho ,\pi ,\sigma ,\nu } \right) \in L^\infty  \left( {Q_t} \right)^2 \times L^\infty  \left( {S_t} \right)^2$ is locally Lipschitz continuous.


}

{\sc {\bf Proof}}. \

To prove this theorem we observe that the IBVP \eqref{ex(2-1-equaz-in-rho)}-\eqref{ex(2-1-condiz-nu-j)} without integral terms is a particular case of \eqref{ex(1-1-IBVP-eq-y)}-\eqref{ex(1-1-dato-frontiera-w)}; therefore we should only study the regularity of the integral terms that appear in the equations \eqref{ex(2-1-equaz-in-pi)}-\eqref{ex(2-1-equaz-in-nu)} and obtain some estimates about these integrals. Now, if we assume
\begin{equation}\label{ex(9-ipotesi-densita-bar)}
\left( {\overline \rho ,\overline \pi,\overline \sigma ,\overline \nu } \right) \in L^\infty  \left( {{0,t_1 };W^{1,\infty } \left( {\Omega } \right)} \right)^2 \times L^\infty  \left( {{0,t_1 };W^{1,\infty } \left( {\Omega_M } \right)} \right)^2,
\end{equation}
then it is not difficult to check that
\begin{equation}\label{ex(9-reg-termini-intg)}
\left[ {\overline \pi   - \pi _l (T)} \right]F_l (\overline \sigma ),\left[ {\overline \pi   - \pi _s (T)} \right]F_s (\overline \nu ) \in 
L^\infty (0,t_1 ;L^{\infty } (\Omega ))
\cap L^1 (0,t_1 ;W^{1,\infty } (\Omega )),
\end{equation}
$$
S_a (\overline \pi  ,\overline \sigma  ,\overline \nu  ),S_q (\overline \sigma  ,\overline \nu  ),N_q (\overline \sigma  ,\overline \nu  ) \in 
L^\infty (0,t_1 ;L^{\infty } (\Omega _M ))
\cap L^1 (0,t_1 ;W^{1,\infty } (\Omega _M )).
$$
Now, we have the following useful estimates, for example, about $Q_l(\overline \sigma, \overline \sigma)$
\begin{equation}\label{ex(9-stime-termini-intg)}
\left\| {Q_l (\overline \sigma  ,\overline \sigma  )} \right\|_{L^1 (0,t;L^\infty  (G))}  \le Ct\left\| {\overline \sigma  } \right\|_{L^\infty  (S_t )}^2, 
\end{equation}
$$
\left\| {\nabla _{(m,x)} Q_l (\overline \sigma  ,\overline \sigma  )} \right\|_{L^1 (0,t;L^\infty  (G))}  \le Ct(\left\| {\overline \sigma  } \right\|_{L^\infty  (S_t )}^2  + \left\| {\nabla _{(m,x)} \overline \sigma  } \right\|_{L^\infty  (S_t )}^2 ), \quad 0< t \le t_1.
$$
Moreover, assuming $\overline \sigma^{(k)} \in L^\infty  \left( {{0,t_1 };W^{1,\infty } \left( {\Omega_M } \right)} \right)$ $(k=1,2)$, we deduce 
\begin{equation}\label{ex(9-stime-pertur-termini-intg)}
\left\| {Q_l (\overline \sigma  ^{(2)} ,\overline \sigma  ^{(2)} ) - Q_l (\overline \sigma  ^{(1)} ,\overline \sigma  ^{(1)} )} \right\|_{L^1 (0,t;L^\infty  (G))}  \le  
\end{equation}
$$
\le Ct(\left\| {\overline \sigma  ^{(1)} } \right\|_{L^\infty  (S_t )}  \vee \left\| {\overline \sigma  ^{(2)} } \right\|_{L^\infty  (S_t )} )\left\| {\overline \sigma  ^{(1)}  - \overline \sigma  ^{(2)} } \right\|_{L^\infty  (S_t )}.
$$
In a similar way, we obtain analogous estimates for the other integral terms.

After having obtained these regularity results and estimates for the integral terms, we understand that it is not so hard to extend the proof of the Theorem 3.1 to prove Theorem 8.1.

\  $ \square$

\end{document}